\documentclass{amsart}

\usepackage{epsfig} 
\usepackage{graphics}

\usepackage{amssymb}
\usepackage{amscd}
\usepackage{amsmath}
\usepackage{enumerate}

\newtheorem{theo}{Theorem}
\newtheorem{lema}[theo]{Lemma}
\newtheorem{cor}[theo]{Corollary}
\newtheorem{prop}[theo]{Proposition}

\newtheorem{definition}[theo]{Definition}

\newtheorem{example}[theo]{Example}

\newcommand{\CC}{{\mathbb{C}}}
\newcommand{\NN}{{\mathbb{N}}}

\newcommand{\RR}{{\mathbb{R}}}
\newcommand{\SSS}{{\mathbb{S}}}

\newcommand{\calA}{{\mathcal{A}}}

\newcommand{\comp}{{\circ}}

\begin{document}
\title[Equisingularity at the normalisation]{Equisingularity at the normalisation}
\author{Javier Fern\'andez de Bobadilla}
\address{Departamento de \'Algebra. Facultad de Ciencias Matem\'aticas. Universidad Complutense. Plaza de Ceincias 3. 28040 Madrid. Spain.}
\email{javier@mat.uned.es}
\author{Mar\'ia Pe Pereira}
\address{Departamento de \'Algebra. Facultad de Ciencias Matem\'aicas. Universidad Complutense.  Plaza de Ceincias 3. 28040 Madrid. Spain.}
\email{maria.pe@mat.ucm.es}
\thanks{First author supported by Ramon y Cajal contract. Both authors supported by the Spanish project MTM2004-08080-C02-01.}
\dedicatory{Dedicated to Jos\'e Mar\'ia Montesinos, on the ocassion of his 60th birthday}
\date{}



\maketitle

\section{Introduction}

We look at topological equisingularity of a holomorphic family of reduced mapping germs $f_t:(\CC^3,O)\to \CC$ over a contractible base $T$ having non-isolated singularities, by means of their normalisations. 
The relation between surfaces in $\CC^3$ with non-isolated singularities and normal surface 
singularities via the normalisation is a very convenient way of exchanging information between both categories. This has been already exploited fruitfully by T. de Jong and D. van Straten in order 
to study deformations of singularities~\cite{JS1},~\cite{JS2},~\cite{JS3}.   

Moreover we apply our results to the study of topological $\calA$-equisingularity of parametrised surfaces. The main observation is that the normalisation of a parametrised surface coincides with the parametrisation itself. 

We introduce the notion of Equisingularity at the normalisation for a family $f_t$ (see Definition \ref{equinorm}). It turns out that in many cases, equisingularity at the normalisation characterises topological embedded equisingularity and $R$-equisingularity. More precisely we prove the following theorem:\\

\noindent \textbf{Theorem A.} {\it 
\begin{itemize}
\item If for a single $t\in T$ each of the generic transversal singularities of $f_t$ has at least one smooth branch, then equisingularity at the normalisation is equivalent to topological $R$-equisingularity.
\item If L\^e's Conjecture holds ~\cite{Le}, \cite{Bo2} (see Section 4 for the statement), then equisingularity at the normalisation is equivalent to topological $R$-equisingularity without any assumption on transversal singularities.  
\end{itemize}}

\vspace{.25cm}

The fact that topological $R$-equisingularity implies equisingularity at the normalisation is not difficult and it is proved at the begining of the paper. The proof of the converse splits in two parts. In the first part we prove the following:\\

\noindent\textbf{Theorem B.} {\it
Equisingularity at the normalisation and equisingularity at the singular locus implies topological $R$-equisingularity.}\\

The idea of the proof is based in \cite{Bo}. Here we also use the notion of {\it cuts}, which are real hypersurfaces in $T\times \CC^3$ which enclose a neighbourhood of $T$ and which imitate the situation of having a uniform radius for the Milnor fibration of $f_t$. Due to the absence of h-cobordism Theorem in the dimension we are dealing with, the tecniques needed here are of a different nature. The reader may see equisingularity at the normalisation as a substitute for the h-cobordism Theorem. In the second part we search for hypothesis under which equisingularity at the normalisation implies equisingularity at the singular locus. These hypothesis are those which appeared in Theorem A.   

To each family of parametrisations $\varphi_t:(\CC^2,O)\to (\CC^3,O)$ we can associate a family of functions $f_t:(\CC^3,O)\to \CC$ whose zero set coincide with the images of the family $\varphi_t$. In this situation, equisingularity at the normalisation means the constancy of the Milnor number at the origin of the inverse image by $\varphi_t$ of the singular set of the image of $\varphi_t$. We prove the following theorem:\\

\noindent\textbf{Theorem C.} {\it
Topological $\calA$-equisingularity of $\varphi_t$ implies the constancy of the Milnor number at the origin of the inverse image by $\varphi_t$ of the singular set of the image of $\varphi_t$. Conversely, 
\begin{itemize}
\item  if for a single $t\in T$ each of the generic tranversal singularities of $\varphi_t(\CC^2)$ has at least one smooth branch, then the constancy of the Milnor number at the origin of the inverse image by $\varphi_t$ of the singular set of the image of $\varphi_t$ implies topological $\calA$-equisingularity.
\item If L\^e's Conjecture is true, the conclusion holds in general.  
\end{itemize}}

For the particular case of finitely determined mapping germs (which authomatically satisfy that generic transversal singularities are of Morse type), the result was obtained also in ~\cite{CHR}. This theorem answers a question of D. Mond in~\cite{Mo}. The proof of Calleja-Bedregal, Houston and Ruas uses 
a technique of Gaffney~\cite{Ga} which has the virtue of giving Whitney equisingularity as well. The finite-determinacy of germs is essential 
for the use of that technique. In fact we expect that in the infite-determined case topological $\calA$-equisingularity does not
necessarily imply Whitney equisingularity. Certainly equisingularity at the normalisation does not imply
Whitney equisingularity: Briancon and Speder's family is a topologically equisingular family of normal 
surfaces in $\CC^3$ (and hence equisingular at the normalisation), which is not Whitney equisingular. On the other hand our technique works under much weaker hypothesis because it uses purely topological tools instead of integration of vector fields.

Our methods also help in the study of $\mu$-constant deformations in $\CC^3$ since the normalisation of a $\mu$-constant family is the identity map.
Thus, Theorem~B has the following corollary:\\

\noindent\textbf{Corollary D.} {\it If a $\mu$-constant family of surfaces in $\CC^3$ is topologically abstractly equisingular, then it is topologically $R$-equisingular.}\\

As a final remark, we notice that our topological trivialisations are {\it global} over the original base $T$ as long as the base is contractible and compact.  

\section{Setting}
Let $T$ be a contractible and compact complex manifold with a base point that is denoted by $0\in T$.
Let $U$ be  a neighbourhood of $T$ in $T\times \CC^3$; we denote $T\times \{0\}$ simply by $T$.

Throughout the article we will consider a holomorphic function
$$f:U\to\CC$$
with $T\subset f^{-1}(0)$. As we denote with $f_t$ the restriction of
$f$ to $\{t\}\times\CC^3$, we will speak of a holomorphic family $f_t$ refering to this situation.
We always ask $f$ and $f_t$ be reduced germs 
at the origin for any $t\in T$.

We will call $\tau:T\times\CC^3\to T$ the projection onto $T$  and denote $$F:=(\tau,f):U\to T\times \CC.$$
Let $(z_1,...,z_n)$ be a coordinate system of $\CC^n$. Define $C:=V(\partial f/\partial z_1,...,\partial f/\partial z_n)$.
We call $X:=f^{-1}(0)$ the zero set of $f$ and $\Sigma:=C\cap X$. Analogously we call $X_t:=f_t^{-1}(0)$ and $\Sigma_t:=\mathrm{Sing} X_t$. Then $\Sigma$ is equal 
 to the union of the $\Sigma_t$.

In general, given any subset $A$ contained in $U$ we will denote $A\cap\tau^{-1}(T')$ by $A_{T'}$ for any subset $T'$ of $T$. In particular, if $T'=\{t\}$ we will denote it by $A_t$.\\
For any function $g:T\times A\to B$, we will denote the function $g(t,\_):A\to B$ by $g_t$.

When we are working in some $\{t\}\times \CC^3$ and the context is clear enough we will simply figure out working in $\CC^3$.  

We consider the following equisingularity notions:

\begin{definition}
\label{equiclassic}
The family $f_t$ is {\em topologically equisingular} if there is a homeomorphism germ
\begin{equation}\label{topequi}\Phi:T\times(\CC^3,O)\to T\times(\CC^3,O),\end{equation}
such that we have the equality $\tau=\tau\comp \Phi$ and $\Phi$ sends $X$ to $T\times X_0$. 
 If in addition $\Phi$ preserves fibres of $f$, we say that $f_t$ is {\em topologically $R$-equisingular} (the prefix $R$ stands for right equivalence).
If moreover, the restrictions $\Phi|_{T\times \CC^3\setminus \Sigma}$, $\Phi_{\Sigma\setminus T}$ are smooth, we say that $f_t$ is {\em piecewise smoothly $R$-equisingular}.

The family $f_t$ is {\em abstractly topologically equisingular} if there is a homeomorphism germ
\[\Psi:(X,T)\to T\times(X_0,O)\]
such that we have the equality $\tau|_X=\tau\comp \Psi$. If in addition, $\Psi$ sends $\Sigma$ to $T\times \Sigma_0$, we say that $f_t$ is {\em abstractly topologically equisingular in the strong sense}.
\end{definition}
Topological $R$-equisingularity is stronger than topological equisingularity, which is stronger than strong abstract
topological equisingularity, which is stronger than abstract topological equisingularity.
\sffamily
\section{Equisingularity at the normalisation}
We will relate the previous equisingularity notions with the following; consider the normalisation mapping 

\[n:\hat X\to X.\]

Consider the set $\hat{\Sigma}:=n^{-1}(\Sigma)$ in which the normalisation mapping is not a local isomorphism.

In general, for any subspace $A$ of $X$, the space $n^{-1}(A)$ is denoted by $\hat A$ and for any function $g:X\to B$, the composition $\hat g:=g\comp n:\hat X\to B$ is denoted by $\hat g$.

\begin{definition}
\label{equinorm}
The family $f_t$ is {\em equisingular at the normalisation} if there is a homeomorphism
$$\alpha:(\hat X,\hat \Sigma,\hat T)\to T\times(\hat X_0,\hat \Sigma_0,n_0^{-1}(O))$$
such that we have the equality $\hat \tau=p_1\comp\alpha$ where $p_1$ is the projection of  $T\times \hat X_0$ to the first factor.
If in addition, the restrictions $\alpha|_{\hat X\setminus \hat\Sigma}$, $\alpha|_{\hat\Sigma\setminus n^{-1}(O)}$ are smooth
we say that $f_t$ is {\em piecewise smoothly equisingular at the normalisation}.
\end{definition}
(To avoid confussion we do not use the notation $\hat O$ or $\hat O_t$.)

\noindent\textbf{Notation.} By equisingularity at the normalisation it is inmediate that, if we index the irreducible components of $\hat\Sigma$ as $\hat\Sigma^1$,...,$\hat\Sigma^r$, then the components of $\hat\Sigma_t$ can also be indexed as $\hat\Sigma_t^1$,...,$\hat \Sigma_t^r$ paralelly in the sense that $\hat\Sigma_t^i$ is contained in $\hat\Sigma^i$ for any $t\in T$.

Here are some consequences of equisingularity at the normalisation 

\begin{lema}
\label{consecnorm}
If $f_t$ is equisingular at the normalisation then
\begin{enumerate}[(i)]
\item the number of irreducible components of the germ $(X,(t,O))$ is independent of $t$.
\item The mapping $n_t:\hat X_t\to X_t$ is the topological normalisation, that is, the continuous mapping underlying it,
coincides with the continuous mapping underlying the normalisation of $X_t$.

\item Each connected component of $\hat X_t$ contains a unique point of $n_t^{-1}(O)$. 

\item For any $t\in T$, the space $(\hat X_t)_{red}$ is smooth outside $n_t^{-1}(O)$.
\item For any $t\in T$, the curve $(\hat\Sigma_t)_{red}$ is smooth outside $n_t^{-1}(O)$.
\item If $\hat\Sigma_t$ is a flat, reduced family of curves then $(\hat \Sigma_t,n_t^{-1}(O))$ is a $\mu$-constant family of multi-germs. 

\item If $\hat X_t$ is equal to $\CC^2$ for any $t$, then $\hat\Sigma_t$ is a $\mu$-constant family.
\end{enumerate}
\end{lema}

\begin{proof}
 The number of irreducible components of the germ $(X,(t,O))$ coincides with the number of preimages of $(t,O)$ by the normalisation mapping, which is the cardinality of $n^{-1}(t,O)$. As $n^{-1}(t,O)$ is equal to $\alpha_t(n^{-1}(0,O))$, the first claim holds.

 Since $n_t:\hat X_t\to X_t$ is finite and generically $1:1$, due to the uniqueness of the topological normalisation, 
it is enough to check that $\hat X_t$ is topologically normal at any point.
 By definition, the surface $\hat X_t$ is topologically normal if and only if $\hat X_t\setminus \mathrm{Sing}(\hat X_t)$ is locally connected at each point of $\mathrm{Sing}(\hat X_t)$.
Suppose that this is not the case. The set $A$ of points of $\mathrm{Sing} \hat X_t$ where this fails, is the 
 locally closed analyitic subset of points of $\hat X_t$ where $\hat X_t$ is not locally irreducible 
By equisingularity at the normalisation we can take a representative $\hat X$ homeomorphic through $\alpha$ to the product $T\times \hat X_t$. As any such homeomorphism transforms $T\times A$ into $\mathrm{Sing}(\hat X)$, 
we obtain that $\hat X\setminus\mathrm{Sing}(\hat X)$ is not locally connected at $\mathrm{Sing}(\hat X)$. This contradicts the normality of $\hat X$ and proves $(ii)$.

 Choosing a small enough representative of $X$, we can assume that $(iii)$ is true 
for $t=0$. For $t\neq 0$ it follows from the existence of $\alpha_t$.

 As $\hat X_0$ is the topological normalisation of $X_0$, it has only isolated singularities. As $\hat\Sigma_0$ is a curve, it has only isolated
singularities. Choosing a sufficiently small representative of $X$ we may assume that $\mathrm{Sing}(\hat X_0)\cup \mathrm{Sing}(\hat\Sigma_0)$ is included 
in $n_0^{-1}(O)$. Now we use the existence of the homeomorphism
$\alpha_t$: by Mumford's Theorem~\cite{Mu}
 the surface $\hat X_t$ is smooth outside $n_0^{-1}(O)$; anologously, the curve $\hat \Sigma_t$ can not be singular
outside $n_0^{-1}(O)$, because we would have a plane curve singularity germ topologically equivalent to a smooth one.
This proves $(iv)$ and $(v)$.

The last assertion is a direct consequence of the sixth, which
follows by the definition of equisingularity at the normalisation and
the characterisation of equisingular families of flat, reduced families of curves given in~\cite{BuG}. 
\end{proof}

\noindent\textbf{Notation.}
Observe that the image of each component $\hat \Sigma_t^i$ by $n$ is  
an irreducible component of $\Sigma_t$. The restriction
\[n|_{\hat \Sigma^i}:\hat \Sigma^i\to n(\hat \Sigma^i)\]
is a finite mapping 
Denote by 
\[\varphi^i:\widetilde{n(\hat \Sigma^i)}\to n(\hat \Sigma^i)\]
the normalisation mapping. As $\hat \Sigma^i$ is topologically normal (and analytically normal outside $\hat T$) there is a unique finite continuous (analytic outside $\hat T$) mapping
\[\tilde{n}^i:\hat \Sigma^i\to \widetilde{n(\hat \Sigma^i)}\]
such that $n|_{\hat \Sigma^i}=\comp\varphi^i\comp\tilde{n}^i$. 
The mapping $\varphi^i$ is generically $1:1$. 
We denote by $d_{t}^i$ the degree of $\tilde{n}^i_t$. 
We define an equivalence relation $R_t$ in the set $\{1,...,r\}$
defining $i,j$ to be equivalent if we have the equality $n_t(\hat \Sigma_{t}^i)=n_t(\hat \Sigma_{t}^j)$. Then we prove the following lemma
 
\begin{lema}
\label{lemamaria}
If $f_t$ is equisingular at the normalisation then 
\begin{itemize}
\item the relation $R_t$ is independent of $t$, and thus we can speak of a unique relation $R$.
\item For any $t\in T$ and for any $i$ the mapping $\varphi^i_t$ is generically $1:1$.
\item The degrees $d_{t}^i$ are independent of $t$ and thus, we can denote them with $d^i$.
\item The fibres of $\tau|_\Sigma:\Sigma\to T$ are generically reduced.
\end{itemize}
\end{lema}

\begin{proof}
Choose a representative $X$ small enough such that 
\begin{enumerate}[(a)]
\item the singular locus $\Sigma_0$ is a union of topological discs only meeting at the origin $O$,
\item the transversal Milnor number at any point of $x\in \Sigma_0\setminus\{O\}$ of the surface $X_0$ only depends on the connected
component of $\Sigma_0\setminus\{O\}$ to which $x$ belongs, 

\item for any irreducible component $\hat\Sigma_0^i$, the mapping $n|_{\hat \Sigma^i_0}:\hat \Sigma^i _0\to \Sigma_0$ is a ramified covering over its image 
which only ramifies at $n^{-1}(O)$.
\end{enumerate}

The locus of $T$ where $i\sim_{R_t} j$ is the set of points $t$ such that $n(\hat\Sigma^i_t)=n(\hat\Sigma^j_t)$. Define $\Sigma_i:=n(\hat\Sigma_i)$ for any $i$. Consider the restriction 
\[\tau|_{\Sigma^i\cap\Sigma^j}:\Sigma^i\cap\Sigma^j\to T.\]
As $\Sigma^i$ and $\Sigma^j$ are irreducible, if the dimension of the fibre $(\Sigma^i\cap\Sigma^j)_t$ 
is not zero then $\Sigma^i_t$ is equal to $\Sigma^j_t$ and the fibre coincides with them. 
By the upper semicontinuity of the fibres of a morphism the locus where $i\sim_{R_t} j$ is closed 
analytic in $T$.

Let $A\subset n(\hat\Sigma^i)$ be the closed analytic subset where $n(\hat\Sigma^i)$ is not topologically
normal. The locus of $T$ where $\varphi^i_t$ is not generically $1:1$ is the set where the dimension of 
the fibres of $\tau|_A:A\to T$ jumps to $1$. This is a proper closed analytic subset.  

Let $B\subset\widetilde{n(\hat\Sigma^i)}$ be the branching locus of the branched cover
$\tilde{n}^i:\hat \Sigma_{t}^i\to \widetilde{n(\hat \Sigma_{t}^i)}$. The degree $d^i_t$ is lower-semicontinuous in $t$, and the locus where it jumps down is closed analytic because it is the locus where the 
dimension of the fibres of the mapping $\tau|_B:B\to T$ jumps up to $1$. 

We have proved that there is a proper closed analytic subset $\Delta\subset T$ such that $R_t$ and  $d^i_t$ for $1\leq i\leq r$ are constant and $\varphi^i_t$ is generically $1:1$ outside $\Delta$. Moreover
$R_t$ is finer that $R_{t'}$ and $d^i_{t}\geq d^i_{t'}$ if $t'\in\Delta$ and $t\in T\setminus\Delta$.
In particular $R_t$ and $d^i_t$ are generically constant.

It is clearly enough to treat the case when $T$ is a disc and $0$ is the unique point over which the 
statement of the Lemma might be false.

\noindent{\textbf Claim}. Let $l:\CC^3\to\CC$ be a generic linear function vanishing at $O$. Denote with $L_s$ the hyperplane $l=s$. By genericity,
$L_s$ meets $\Sigma_0$ transversely for any $s\neq 0$ small. Then there exists a neighbourhood $T'$ of $0$ in $T$ such that
\begin{equation}
\label{familia}
\tau:(X\cap (T'\times L_s),\Sigma\cap (T'\times L_s))\to T'.
\end{equation}
 is a $\mu$-constant family of multigerms. 
 
Let us prove the claim. As $n_0|_{\hat \Sigma_0}:\hat\Sigma_0\to\Sigma_0$ is an unramified covering outside
$n_0^{-1}(O)$ and $\Sigma_0$ meets $L_s$ transversely, the mapping $n_0|_{\hat\Sigma_0\setminus\{O\}}$ 
is transverse to $L_s$; thus $n_0:\hat X_0\to\CC^3$ is transverse to $L_s$ at any point of 
$\hat\Sigma_0\setminus\{O\}$. Consequently 
all the components of the multi-germ
$(n_0^{-1}(L_s),\hat \Sigma_0\cap n_0^{-1}(L_s))$ are smooth germs of curves. From now on we fix a generic $s$ with this property.

By equisingularity at the normalisation the family $\tau:\hat \Sigma\to T$ is topologically equisingular. This implies easily that any fibre
$\hat \Sigma_t$ is reduced at its generic point. 
As a consequence we find that the mapping
\[\hat\tau:\hat\Sigma\cap n^{-1}(T\times L_s)\to T\]
is a local diffeomorphism near $\hat\Sigma_0\cap n_0^{-1}(L_s)$. Thus the mapping
\[\hat\tau:(n^{-1}(T\times L_s),\hat \Sigma\cap n^{-1}(T\times L_s))\to T\]
is a submersion at a neighbourhood $U$ of $\hat \Sigma_0\cap n_0^{-1}(L_s)$ in $n^{-1}(T\times L_s)$.
Taking local charts for the submersion we can choose $U$ and a neighbourhood $T'$ of $0$ in $T$ such 
that 
\[\hat\tau:(U,\hat \Sigma\cap U)\to T'\]
is a locally trivial fibration whose fibre is diffeomorphic to the disjoint union of several puctured discs $(D,*)$. 
The image $(V,\Sigma\cap (T'\times L_s)):=n(U, \hat \Sigma\cap U)$ is a representative of the family of multi-germs of plane curve singularities
$$\tau:(X\cap (T'\times L_s)),\Sigma\cap (T'\times L_s))\to T'.$$
The mapping
\[n|_V:(U,\hat \Sigma\cap U)\to (V,\Sigma\cap (T'\times L_s))\]
is a simultaneous (non-embedded) resolution of singularities of the family over $T'$. 
Using this fact, we are going to prove that
the Milnor number of the family of multi-germs is constant.

Observe that $\mathrm{Sing}(V_t)$ is equal to $\Sigma_t\cap L_s$ for any $t\in T'$.
The irreducible components of the multigerm $(V_t,\Sigma_t\cap L_s)$ are in a 1-1 correspondence with the discs of $U_t$ and then its cardinality is constant for any $t$. We give an indexing to these components so that for any index $i$ and any $t,t'\in T'$ the components $V_t^i$ and $V_{t'}^i$ are 
images of discs lying in the same connected component of $U$.

Next, we prove
 that the intersection multiplicity in $L_s$ of any pair of corresponding components at $\mathrm{Sing}(V_t)$ is independent of $t$.
If the contrary holds, there exist a pair $(V_t^1,p_t^1)$, $(V_t^2,p_t^2)$ of components of the multigerm $V_t$ such that the following strict inequality holds 
 $$I_{p_t^1}(V_t^1,V_t^2)<I_{p_0^1}(V_0^1,V_0^2)$$  for $t\neq 0$. 
By the conservation of the intersection number 
, the set $V_t^1\cap V_t^2$ contains of a non-empty finite set of points $\{q_t^1,...,q_t^r\}$ different
from $p_t^1$ and convergent to $p_0^1$ as $t$ approaches $0$. Each point $q_t^i$ is a singularity 
 of the multigerm $V_t$. As the singular locus of 
$V_t$ is equal to $\Sigma_t\cap L_s$, and in a small enough neighbourhood $W$ of $p_0^1$ in $V$, the only point of $\Sigma_t\cap L_s\cap W$ is $p_t^1$ 
we have got a contradiction.

Lastly, let us see that the Milnor numbers of the irreducible components are constant. Suppose that there exists a component $(V_t^i,p_t)$ such that we have the strict inequality $\mu(V_t^i,p_t)<\mu(V_0^i,p_0)$
for $t\neq 0$. Then 
$\delta(V_t^i,p_t)<\delta(V_0^i,p_0)$. As $V_t$ is a family of images of parametrisations, the sum of $\delta$-invariants must be preserved. Therefore $V_t$ must have singular points outside $p_t$ which converges to $p_0$ as $t$ approaches $0$. So we get a contradiction and finish the proof of the claim.
 
Now, for any $t\in T'$ we define an equivalence relation $R'_t$ in the finite set $\hat \Sigma_t\cap U$ imposing that two points are equivalents if they
have the same image by $n_t$. The fact that the relations $R'_t$ are independent of $t$ follows from the equisingularity of
the family (\ref{familia}) of plane curves. 
This implies easily 
 that the fibres of $\tau|_\Sigma\to T$ are reduced at their generic points, that $R_t$ is independent
of $t$, and that $\varphi^i_t$ is generically $1:1$ for any $t\in T$. Given that the relation $R_t$ is independent of $t$ and that $\varphi^i_t$ is generically $1:1$ for any $t\in T$, the cardinality of the 
equivalence classes of $R'_t$ only depends on the degrees $d^i_t$. 
Thus $d^i_t$ is independent of $t$ for
any $1\leq i\leq r$. 
\end{proof}

\begin{lema}
\label{lift}
Let $h:X_1\to X_2$ be a homeomorphism between reduced analytic spaces preserving the singular locus.
Let $n_1:\hat X_1\to X_1$ and $n_2:\hat X_2\to X_2$ be the normalisation mappings.
There is a unique homeomorphism $g:\hat X_1\to \hat X_2$ lifting $h$.
\end{lema}
\begin{proof}
Since
\[n_i|_{\hat X_i\setminus \hat \Sigma_i}:\hat X_i\setminus \hat\Sigma_i\to X_i\setminus \Sigma_i\]
is a homeomorphism for $i=1,2$  
We are forced to define
\[g|_{\hat X_1\setminus \hat \Sigma_1}:=n_2|_{\hat X_2\setminus \hat\Sigma_2}^{-1}\comp h|_{X_1\setminus \Sigma_1}\comp n_1|_{\hat X_1\setminus \hat\Sigma_1}.\]

We only have to show that $g|_{\hat X_1\setminus \hat\Sigma_1}$ admits a unique continuous extension to $\hat X_1$.
Uniqueness is clear since $\hat X_1\setminus \hat\Sigma_1$ is dense in $\hat X_1$. Consider $y_1\in \hat\Sigma_1$. Define
$x_1:=n_1(y_1)$ and $x_2:=h(x_1)$. We have to define $g(y_1)$ within the finite set $n_2^{-1}(x_2)$. The set
$n_i^{-1}(x_i)$ is in a $1:1$ correspondence with the irreducible components of the germ $(X_i,x_i)$ for $i=1,2$. Thus there is a
unique irreducible component $(Z,x_1)$ of $(X_1,x_1)$ corresponding to $y_1$. We define $g(y_1)=y_2$, being $y_2$ the unique point
of $n_2^{-1}(x_2)$ corresponding to $(h(Z),x_2)$. It is easy to check that this choice gives a continuous extension.
\end{proof}

\begin{prop}
\label{abstonorm}
If $f_t$ is abstractly topologically equisingular in the strong sense then $f_t$ is equisingular at the normalisation.
\end{prop}
\begin{proof}
If $\Psi$ is the homeomorphism which gives the abstractly topological equisingularity, by Lemma~\ref{lift} there is a unique homeomorphism
\[\alpha:T\times \hat X_0\to \hat X\]
lifting $\Psi$.
\end{proof}

\section{Equisingularity at the singular locus}\label{ssequiatsing}
\begin{definition}
\label{equiatsing}

We say that $f_t$ is {\em equisingular at the singular locus} if
\begin{enumerate}[(i)]
\item the projection $\tau:\Sigma\to T$ is a (non-necessarily reduced) topologically equisingular family of curves.
\item for any $t\in T$, the transversal Milnor number of $f_t$ at {\em any} point of $\Sigma_t\setminus\{O\}$ only depends on the
connected component of $\Sigma\setminus T$ in which it lies.
\end{enumerate}
\end{definition}

Observe that by the first condition of the definition,
the inclussion defines a one to one correspondence between the connected components of
$\Sigma_t\setminus\{O\}$ and the connected components of $\Sigma\setminus T$ for any $t\in T$. 

Example~(51) of~\cite{Bo} shows that the family $\Sigma_t$ is not necessarily reduced 
at the origin.

In equisingularity at the singular locus we only care about the behaviour of the family of functions at the singular locus $\Sigma$ of the central fibre $X$  minus the origin. We will need also to control it at the whole critical locus $C$ minus the origin. Thus, following \cite{Bo}, we make the following definition:

\begin{definition}
\label{equiatcrit}
We say that $f_t$ is {\em equisingular at the critical locus} if
\begin{enumerate}[(i)]
\item it is equisingular at the singular locus and  
\item the spaces $\Sigma$ and $C$ coincide.

\end{enumerate}
\end{definition}

In Section \ref{ssequiatcrit} we will prove that assuming equisingularity at the normalisation, equisingularity at the singular locus implies equisingularity at the critical locus. 

We recall L\^e's conjecture~\cite{Le} (see~\cite{Bo2} for equivalent formulations):\\

\noindent\textbf{Conjecture}. (L\^e)
{\em The image of an injective analytic mapping germ from $(\CC^2,O)$ to $(\CC^3,O)$ is a surface with smooth
singular locus and constant transversal Milnor number.}\\

\begin{prop}
\label{equinormsing}
The following statements hold:
\begin{enumerate}[(i)]
\item If L\^e's conjecture holds, equisingularity at the normalisation implies equisingularity at the singular locus.
\item If for a single $t\in T$ each of the generic transversal singularities of $f_t$ has at least one smooth branch, then the conclusion holds independently of L\^e's conjecture.
\end{enumerate}
\end{prop}

Let us state separately the following important special case of the second assertion:

\begin{cor}
\label{morsetype}
If for a single $t\in T$ the generic transversal singularities of $f_t$ are of Morse type, then equisingularity at the
normalisation implies equisingularity at the singular locus.
\end{cor}

\begin{proof}[Proof of Proposition~\ref{equinormsing}]\quad

\noindent\textbf{Step 1}. We prove that the projection $\tau|_{\Sigma}:\Sigma\to T$ is a topologically equisingular family of curves,
which are smooth outside of the origin. By Lemma~\ref{lemamaria} the fibres of
$\tau:\Sigma\to T$ are generically reduced.

By conditions $(a)-(c)$ of the proof of Lemma~\ref{lemamaria} and equisingularity at the normalisation we can suppose that $\hat\Sigma_t^i$ are discs for any $t$.

As $\hat\Sigma_t^i$ is normal we can factorise the finite mapping \[n|_{\hat\Sigma_t^i}:\hat\Sigma_t^i\to n(\hat\Sigma_t^i)\] through the normalisation $\widetilde{n(\hat\Sigma_t^i)}\to n(\hat\Sigma_t^i)$. 
We obtain a branched covering
\[\hat\Sigma_t^i\to \widetilde{n(\hat\Sigma_t^i)}\]
of degree $d^i$.
Observe that, since $\widetilde{n(\hat\Sigma_t^i)}$ is a connected Riemann surface with boundary, its Euler characteristic must be at most $1$,
reaching the maximal value only when $\widetilde{n(\hat\Sigma_t^i)}$ is a topological disc. An Euler characteristic computation shows
\[1=\chi(\hat\Sigma_t^i)=d_t^i\chi(\widetilde{n(\hat\Sigma_t^i)})-R\]
where $R$ is the total number of ramification points. We deduce that $\widetilde{n(\hat\Sigma_t^i)}$ is a topological disc.

Observe that $\Sigma_t$ consists precisely of the union of the curves $n(\hat\Sigma_t^i)$ and that we know that for any $t\in T'$ all these
discs meet at the origin $O$. For $t=0$ we know that the origin is the unique intersection point of any of these two discs, and that
any of them has no self intersections. In
order to conclude the proof of Step 2 we only have to show that the same holds for any $t\in T$. Suppose the contrary: either
two of the curves $n(\hat\Sigma_t^i)$ meet or one of them has a self intersection for a certain $t\neq 0$.

We analise first the case in which the intersection point is not the origin: suppose that there are indexes $i$, $j$, not necessarily
distinct, points $x\in \hat\Sigma_t^i$, $y\in \hat\Sigma_t^j$, and neighbourhoods $U_x$ of $x$ in $\hat\Sigma_{t,i}$, $U_y$ of $y$ in $\hat\Sigma_{t,j}$
such that we have the equality $n(x)=n(y)$, and the intersection $n(U_x\setminus\{x\})\cap n(U_y\setminus\{y\})$ is empty.
Consider neighbourhoods $V_x$ of $x$ in $\hat X_t$ and $V_y$ of $y$ in $\hat X_t$. The surfaces $n(V_x)$ and $n(V_y)$ in $\CC^3$ share the 
common point $n(x)=n(y)$. Therefore their intersection $Z:=n(V_x)\cap n(V_y)$ is a curve which is obviously contained in $\Sigma_t$.
Hence $n^{-1}(Z)$ is contained in $\hat\Sigma$. Observe that $x$ belongs to $n^{-1}(Z)$, but that the unique common point of $n^{-1}(Z)\cap V_x$
and $\hat\Sigma_t^i\cap V_x$ is the point $x$. Thus $\hat\Sigma$ is not locally irreducible at $x$ and hence, it is singular. This contradicts
Assertion (iv) of Lemma~\ref{consecnorm}.

Suppose now that we have a component $\hat\Sigma_{t,i}$ having the origin as a self intersection point. Then there exists two distinct points
$x,y\in \hat\Sigma_{t,i}$ such that $n(x)=n(y)=0$. As $\hat\Sigma_{t,i}$ is connected, both points are in the same connected component of $\hat X_t$. As by equisingularity at the normalisation each of the connected components of $\hat X_t$ contains a unique point of $n_t^{-1}(O)$, we get a contradiction.\\

\noindent\textbf{Step 2}. We prove the second condition of equisingularity at the singular locus.
As the set $n|_{\hat\Sigma_t^i}^{-1}(0)$ consists of a single point, we conclude that $n|_{\hat\Sigma_t^i}$ is the $d^i$-th cyclic covering of the disc onto itself totally ramified at $n|_{\hat\Sigma_t^i}^{-1}(0)$.
Thus, $n_t$ is locally injective at any point $y\in \hat X_t\setminus n_t^{-1}(0)$. By L\^e's Conjecture
we conclude that $\Sigma_t$ (which is equal to  $\cup_{i=1}^r n(\hat\Sigma_t^i)$) is smooth outside the origin and that for any neighburhood $W$ of any point of
$\hat\Sigma_t^i\setminus n_t^{-1}(0)$ in $\hat X_t\setminus n_t^{-1}(0)$, the image $n(W)$ is a surface possibly singular only at $n(W\cap \hat\Sigma_t^i)$ with constant transversal Milnor number.

In the proof of Lemma~\ref{lemamaria} we have obtained that, for any $t\in T$, the transversal Milnor number of $f_t$ at a {\em generic} point of
$\Sigma_t\setminus\{O\}$ only depends on the connected component of $\Sigma\setminus T$ in which it lies 
(this is the content of the Claim included in the proof of the lemma). Therefore
it is sufficient to show that for any $t\in T$, the transversal Milnor number of $f_t$ at {\em any} point of $\Sigma_t\setminus\{O\}$
only depends on the connected component of $\Sigma_t\setminus\{O\}$ in which it lies.

Consider any point $x\in\Sigma_t\setminus\{O\}$. The number of irreducible components of the transversal singularity coincides with
the number of preimages of $x$ by $n$, and therefore only depends on the connected component of $\Sigma_t\setminus\{O\}$ to which
$x$ belongs. We have obtained as a consequence of L\^e's Conjecture that the Milnor number of each of these irreducible components is
locally constant, and hence depends only on the connected component of $\Sigma_t\setminus\{O\}$ to which
$x$ belongs. It only remains to be shown that the intersection multiplicities between the different branches of the transversal
singularities are locally constant in $\Sigma_t\setminus\{O\}$. Such intersection multiplicities are upper semicontinuous.

Suppose that there is a special point $x\in\Sigma_t\setminus\{O\}$ in which the intersection multiplicity between two branches jumps.
More precisely, given a linear function $l:\CC^3\to\CC$ meeting $\Sigma_t$ transversely at $x$, the branches of the plane curve
singularity given by the section defined by $l$ at a point $x'\in\Sigma_t$ close to $x$ are in bijective correspondance with
the set of preimages of $x'$ by $n_t$. There are two points $y_1,y_2\in n_t^{-1}(x)$ such that the intersection multiplicity of the
two associated branches of the hyperplane section is strictly bigger than the intersection multiplicity of the corresponding branches
at neighbouring hyperplane sections. For $i=1,2$ we choose a small neighbourhood $V_{y_i}$ of $y_i$ in $\hat X_t$, with a biholomorphim
$\varphi_i:V_{y_i}\to D\times D$ (denoting $D$ the disc) such that
the intersection $\varphi_i(\hat\Sigma_t\cap V_{y_i})$ is equal to $\{0\}\times D$ and that $\varphi_i^{-1}(D\times\{a\})$ for $a\in D$ are the
inverse images by $n$ of the hyperplane sections by $L$. By the conservation of the intersection multiplicity, for a certain $a\in D$, the two branches
$n(\varphi_i^{-1}(D\times\{a\})$ for $i=1,2$ must intersect outside $\Sigma$. Their intersection points are singular points of
$X_t$, and therefore must belong to $\Sigma$. We conclude that the intersection multiplicity is locally constant. This concludes
Step 3 and the proof of Assertion $(i)$.\\

We prove Assertion $(ii)$. In Step 1 we have proved that for any $t\in T$, the transversal Milnor number of $f_t$ at a {\em generic} point
of $\Sigma_t\setminus\{O\}$ only depends on the connected component of $\Sigma\setminus T$ in which it lies.
Using this and the fact that for a single $t$ each of the generic tranversal singularities of $f_t$ at $\Sigma_t\setminus\{O\}$ has at least one smooth branch, we deduce that for {\it any} $t\in T$ every generic tranversal singularity of $f_t$ at $\Sigma_t\setminus\{O\}$ has at least one smooth branch.

We observe that we have used L\^e's conjecture only at Step 2. Using the existence of smooth branches, we can replace the use of L\^e's conjecture by Mumford's Theorem (see \cite{Mu}) which states that the unique normal surface germ having a link
with trivial fundamental group is the smooth germ, and prove that $\Sigma_t$ is smooth outside of the 
origin. Then we can use the fact that L\^e's conjecture holds if the image contains a smooth 
branch~\cite{Bo2} and proceed as in Step~2.
\end{proof}

\section{Topological equisingularity}

The purpose of this section is to prove the following theorem:
\begin{theo}\label{equitop}
Let $f_t$ be equisingular at the singular locus and equisingular at the normalisation, then it is $R$-topologically equisingular.
\end{theo}
To prove this, we follow the idea of \cite{Bo} and base our proof on constructing a family of hypersurfaces in $T\times \CC^3$, that we will call \textit{cuts}, which enclose neighbourhoods of $T$ where the restriction of $\tau$ is a trivial fibration and where, near $\Sigma$, we recognise a $\mu$-constant family of isolated singularities. We use the global version of the Theorem of L\^e-Ramanujan in \cite{Bo} to trivialise this $\mu$-constant family and use the equisingularity at the normalisation to trivialise outside that neighbourhood of $\Sigma$.


For the rest of the section the family $f_t$ will be equisingular at the singular locus and equisingular at the normalisation.

\subsection{Cuts}$\quad$

\noindent\textbf{Notation}. Let $E\to N$ be an arbitrary vector bundle with a metric $\rho$ in its fibres. Given any function $\nu:N\to\RR_+$, for any $N'\subset N$ we call
$$E(N',\nu):=\{(x,v)\in E_{N'}:\rho_x(0,v)< \nu(x)\}.$$
For the bundle $\tau$ 
and its restrictions we use the notation $B(T,\nu)$. 

Assuming equisingularity at the singular locus, we can fix a complex normal bundle $\pi:E\to \Sigma\setminus T$ in $T\times\CC^3$ and a positive function $\mu:\Sigma\setminus T\longrightarrow\RR_+$, such that $E(\Sigma\setminus T,\mu)$ can be embedded in $T\times\CC^3$. The restriction
\[\pi|_{E(\Sigma\setminus T,\mu)}:E(\Sigma\setminus T,\mu)\to \Sigma\setminus T\]
is a fibration with fibre a ball in $\CC^2$.
 Observe that the function $f$ restricted to this bundle is a $\mu$-constant family of curves.

\begin{definition}\label{Xcut}
A topological hypersurface $H$ in $T\times\CC^3$ is a {\it cut} for $\Sigma$ of amplitude $\delta$ if

\begin{enumerate}[(i)]
\item the intersection $H\cap E(\Sigma\setminus T,\mu)$ coincides with $E(H\cap\Sigma,\mu)$,
\item the hypersurface $H$ is topologically transversal to $F^{-1}(t,s)$ outside $\Sigma$, for any $(t,s)\in T\times D_\delta$, 
\item the hypersurface $H$ separates $F^{-1}(T\times D_\delta)$ in several connected components, one of them is a bounded neighbourhood of $T$; we call it $\mathrm{Int} H(\delta)$. (We will not specify the amplitude if it is clear from the context). Besides, $\overline{\mathrm{Int} H(\delta )}$ is a topological manifold with boundary such that $$\tau|_{\overline{\mathrm{Int} H(\delta)}}:\overline{\mathrm{Int} H(\delta)}\to T$$ is a topological submersion.
\item for any $t\in T$, 
$(\mathrm{Int} H(\delta)\cap\Sigma)_t$ is homeomorphic to a cone over $H\cap\Sigma_t$.
\end{enumerate}
\end{definition}

Abusing notation slightly, we consider the subspace $\hat{E}(\mu):=n^{-1}(E(\Sigma\setminus T,\mu))$
which inherits the unique natural bundle structure
\begin{eqnarray}\label{Egorro}\hat{\pi}:\hat{E}(\mu)\to \hat \Sigma\setminus \hat T&&\end{eqnarray}
compatible with the one in $E(\Sigma\setminus T,\mu)$.

\begin{definition}
\label{Ycut}
A topological real hypersurface $\hat{H}$ in $\hat X$ is a {\em cut} for $\hat \Sigma$ in $\hat X$ if
\begin{enumerate}[(i)]
\item The restricition $\hat\tau|_{\hat H}$ is a proper topological submersion.
\item For any $t\in T$ the intersection $(\hat{H}\cap \hat\Sigma)_t$ is a smooth real curve saturated by $n$,
that is $n^{-1}(n((\hat{H}\cap \hat\Sigma)_t))$ coincides with $(\hat{H}\cap \hat\Sigma)_t$.
\item For any $t\in T$ the smooth curve $n(\hat{H}\cap \hat\Sigma)_t$ separates $\Sigma_t$ in several connected components. The only component
containing the origin is a cone over the curve $n(\hat{H}\cap \hat\Sigma)_t$.
\item The intersection $\hat{H}\cap \hat{E}(\mu)$ 
coincides with $\hat{\pi}^{-1}(\hat H\cap \hat\Sigma)$.
\item The cut $\hat{H}$ separates $\hat X$ in several connected components, the union of those meeting $\hat T$ is a neighbourhood of $\hat T$ that we call $ \mathrm{Int}\hat{H}$. Any connected component of $(\mathrm{Int}\hat H)_t$ is a cone over the only point of $n_t^{-1}(0)$ that it contains.
\end{enumerate}
\end{definition}

\noindent\textbf{Note}. The notation of $\hat H$ is due to the fact that, as the reader will see in Section \ref{ssXcut}, a cut $\hat H$ can be the preimage by the normalisation $n$ of a cut $H$.

In order to simplify the exposition we do the following general definition:
\begin{definition}
\label{interiores}
Let $A$ be a subspace of $U$ (resp. $\hat X$) that contains a subset $T'$ of $T$ (resp. $\hat T'\subset \hat T$).
Let $M$ a subspace of $A$ which does not meet $T'$ (resp. $\hat T'$).
We can define
$$(\mathrm{Int}_A M)_{T'}$$ as the union of the connected components of $(A\setminus M)_{T'}$ meeting $T'$ (resp. with $\hat T'$).
$(\mathrm{Int}_A M)_{T'}(\delta_i)$ will stand for $(\mathrm{Int}_{A\cap f^{-1}(D_{\delta_i})} M)_{T'}$.
\end{definition}
When the context makes clear who $A$ and $T'$ are, we will not specify them. Observe that this definition is coherent with the definition of $\mathrm{Int} H$ and $\mathrm{Int} \hat H$ given above for certain cuts $H$ and $\hat H$: $\mathrm{Int} H=(\mathrm{Int}_UH)_T$, $\mathrm{Int} \hat H=(\mathrm{Int}_{\hat X}\hat H)_{\hat{T}}$. Also, as we will do in the Section \ref{ssYcut}, taking a hypersurface $M$ in $\Sigma\setminus T$ we can consider $\mathrm{Int} M=(\mathrm{Int}_\Sigma M)_T$.

\subsection{A glueing technique}
In this subsection we explain a technical tool that we use throughout the rest of the paper. 
In our proofs we will have to construct piecewise several homeomorphisms. We will be working on spaces fibred over the contractible space of parameters $T$ and most of our construtions in each of the pieces will be trivial over the parameter $0\in T$. We will show in this subsection how to use the retraction of $T$ onto $0$ to make them compatible.

Let $(\Lambda,0)$ be a contractible and compact topological space.
Let $A$ be a topological space. 
Let $f_i:A\to B_i$ be continuous mappings to topological spaces
$B_i$, for $i\in\{1,...,m_1\}$ and $C_i$, $D_j$ be subspaces of $A$, for $i\in\{1,...,m_2\}$, $j\in\{1,...,m_3\}$.

Let
\[\tilde h:\Lambda\times A\to \Lambda\times A\]
be a self-homeomorphism of the form $\tilde h:=(Id|_\Lambda,h_\lambda)$  
with $\tilde h_\lambda:A\to A$ a family of homeomorphisms
depending continuously in $\lambda\in\Lambda$. Suppose furthermore that this family of homeomorphisms 
satisfies the following conditions:
\begin{enumerate}[(a)]
\item the homeomorphism $\tilde h_0$ is the identity,
\item the fibres of each of the mappings $f_i:A\to B_i$ are preserved by $\tilde h_\lambda$ for all $\lambda\in\Lambda$, ($f_i\comp \tilde h_\lambda=f_i$),
\item the subspaces $C_i$ are left invariant by $\tilde h_\lambda$ for all $\lambda\in \Lambda$,
\item the subspaces $D_j$ are fixed by $\tilde h_\lambda$ for all $\lambda\in\Lambda$.
\end{enumerate}

We associate any continuous mapping
\[g:(\Lambda,0)\to (\Lambda,0)\]
with a homeomorphism
\[\tilde g:\Lambda\times A\to\Lambda\times A\]
defined by $\tilde g(\lambda,x):=(\lambda,\tilde h_{g(\lambda)}(x))$; the homeomorphisms 
$\tilde g_\lambda$ obviously satisfy conditions (a)-(d).

\begin{lema}
\label{glueprev}
Let $\tilde h:\Lambda\times A\to \Lambda\times A$ be a homeomorphism as above. There exists an isotopy
\[\phi:\Lambda\times A\times I\to\Lambda\times A\times I\]
between $\phi_0=\tilde h$ and $\phi_1=Id_{\Lambda\times A}$
of the form $$\phi(\lambda,x,t)=(\lambda,\phi_{(\lambda,t)}(x),t)$$
for certain family of homeomophisms $\phi_{(\lambda,t)}:A\to A$ which satisfy 
(a)-(d).
\end{lema}
\begin{proof}
As $\Lambda$ is contractible, there exists a strong deformation retract of $\Lambda$ to the point $0$ $$\xi:\Lambda\times I\to\Lambda$$ which satisfies $\xi(\lambda,0)=\lambda$, $\xi(\lambda,1)=0$ for all $\lambda\in \Lambda$ and $\xi(0,t)=0$ for all $t\in [0,1]$.

We define $\phi_t$ as the homeomorphism $\tilde\xi_t$ associated to $\xi_t$ and $\phi:=(\phi_t,Id|_{[0,1]})$.

It is clear that the mapping $\phi$ is bijective and continuous. Its inverse is continuous because the family
of homeomorphisms $\tilde h^{-1}_\lambda$ depends continuously in $\Lambda$ (for being $\tilde h^{-1}:\Lambda\times A\to \Lambda\times A$ continuous.). It is obvious that $\phi$ satisfies the rest of the requirements.
\end{proof}

Let $E$ be another topological space and
 $\beta:E\to\Lambda$ a continuous mapping.
Consider continuous mappings $f_i:E\to B_i$ to topological spaces $B_i$ for $1\leq i\leq m_1$ and
subspaces $C_i\subset E$ for $1\leq i\leq m_2$.

\begin{definition}
\label{trivrel}
A {\em trivialisation of} $\beta:E\to\Lambda$ {\em relative to} $f_i:E\to B_i$ for $1\leq i\leq m_1$ and
$C_i\subset E$ for $1\leq i\leq m_2$ is a homeomorphism
\[\xi:E\to \Lambda\times E_0\] 
with $E_0=\beta^{-1}(0)$ such that
\begin{itemize}
\item we have the equality $\beta=p_1\comp\xi$, where $p_i$ denotes the projection of $\Lambda\times E_0$ to the $i$-th factor,
\item it satisfies $\xi|_{E_0}=(0,Id|_{E_0})$,
\item we have the equality $f_i=f_i|_{E_0}\comp p_2\comp\xi$ for $1\leq i\leq m_1$,
\item we have the equality $\xi(C_i)=\Lambda\times (E_0\cap C_i)$ for $1\leq i\leq m_2$.
\end{itemize}
\end{definition}

\begin{example}
The definition of equisingularity at the normalisation can be stated saying that $\hat\tau:=\tau\comp n:\hat X \to T$ admits a trivialisation $\alpha:\hat X\to T\times\hat X_0$ relative to $\hat\Sigma$ and $\hat T$.
\end{example}

\begin{example}
If $\beta:E\to \Lambda$ is a proper differentiable submersion, then the Ehresman Fibration Theorem tells that there exist a trivialisation of $\beta$.
\end{example}
\begin{example}
Analogously, if $\beta:E\to\Lambda$ is a proper topological submersion, then there exist also a trivialisation of $\beta$ (see~\cite{KS} page 59).
\end{example}

Assume that $E$ is the union of two closed subspaces $E_1$ and $E_2$.
Suppose that for $i=1,2$ we have trivialisations of $\beta|_{E_i}$
\[\xi_i:E_i\to \Lambda\times (E_i)_0\]
relative to $f_i|_{E_i}$ for $1\leq i\leq m_1$
and to the subsets $C_j\cap E_i$ for $1\leq j\leq m_2$ and $E_1\cap E_2$.
We want to modify $\xi_1$ and $\xi_2$ so that they can be plumbed together to give a trivialisation of $E$ relative to
$f_i$ for $1\leq i\leq m_1$ and $C_i$ for $1\leq i\leq m_2$.

\begin{definition}
\label{collar}
Let $A$ be a subspace of $E$. We say that a neighbourhood $C(A)$ of $E$ is a collar for $A$ in $E$
{\em relative to} $f_i$ for $1\leq i\leq m_1$ and to the subsets $C_i$ for $1\leq i\leq m_2$,
if there is a homeomorphism
\[\theta:A\times [0,1]\to C(A)\]
such that
\begin{itemize}
\item its restriction $\theta_0$ is the inclusion of $A$ in $C(A)$,
\item it satisfies $f_i\comp\theta=f_i\comp p_1$ for $1\leq i\leq m_1$, with $p_1:A\times [0,1]\to A$ the projection onto the first factor,
\item we have the equality $\theta((C_i\cap A)\times [0,1])=C_i\cap C(A)$.
\end{itemize}
\end{definition}
\begin{lema}
\label{glueing}
If $A:=(E_1\cap E_2)_0$ has a collar $C(A)$ in $(E_1)_0$
relative to $f_i|_{(E_1)_0}$ for $1\leq i\leq m_1$ and relative to the subsets $C_i\cap (E_1)_0$ for $1\leq i\leq m_2$,
then the trivialisations $\xi_1$ and $\xi_2$ can be plumbed in order to obtain a trivialisation of $\beta:E\to \Lambda$,
$$\xi:E\to\Lambda\times E_0,$$
relative to $f_i$ for $1\leq i\leq m_1$ and $C_i$ for $1\leq i\leq m_2$ and such that
$\xi|_{E_2}=\xi_2$ and
$\xi|_{E_1\setminus\xi_1^{-1}(\Lambda\times{C(A)})}=\xi_1$.

\end{lema}
\begin{proof}
The homeomorphism
\[h:=\xi_2|_{E_1\cap E_2}\comp \xi_1^{-1}|_{\Lambda\times (E_1\cap E_2)_0}:\Lambda\times (E_1\cap E_2)_0\to \Lambda\times (E_1\cap E_2)_0\]
satisfies the hypothesis of Lemma~\ref{glueprev} so there exists an isotopy
$$\phi:\Lambda\times (E_1\cap E_2)_0\times [0,1]\to\Lambda\times (E_1\cap E_2)_0\times [0,1]$$
with $\phi_0=h$ and $\phi_1=Id_{\Lambda\times A}$.

We call $\theta$ the homeomorphism which gives $C(A)$ the structure of collar of $A$ in $(E_1)_0$ relative to $f_i$ for $1\leq i\leq m_1$ and $C_i$ for $1\leq i\leq m_2$.

We can construct $\xi:E\to \Lambda\times E_0$ piecewise as follows

$$\xi|_{E_2}:=\xi_2,$$ 
$$\xi|_{E_1\setminus\xi_1^{-1}(\Lambda\times C(A))}:=\xi_1,$$
$$\xi|_{\xi_1^{-1}(\Lambda\times{C(A)})}:=(Id_{\Lambda},\theta)\comp\phi\comp(Id_{\Lambda},\theta^{-1})\comp\xi_1.$$

\end{proof}

\subsection{Existence of cuts for $\hat \Sigma$}\label{ssYcut}
In proving the existence of cuts for $\hat \Sigma$ we will use several auxiliary hypersurfaces in $\Sigma$ that we construct thanks to the following lemma:
\begin{lema}
\label{Singcut}
Suppose that $f_t$ is equisingular at the singular locus.
Given any positive function $\eta:T\to \RR_+$. Choose $\epsilon_*>0$ such that $\SSS_{\epsilon}$ is transverse to $\Sigma_0$ for any
$\epsilon<\epsilon_*$. There exists $0<\epsilon<\epsilon_*$ and a smooth hypersurface $M$ of $\Sigma\setminus T$ such that
\begin{enumerate}[(i)]
\item the restriction $\tau|_M$ is a submersion,
\item it separates $\Sigma$ in several connected components, the only component $IntM$ containing $T$ is contained in $B(T,\eta)$,
\item for any $t\in T$ the fibre $M_t$ is a disjoint union of circles, one located in each of the irreducible components of $\Sigma_t$,
\item for any $t\in T$ the fibre $IntM_t$ is homeomorphic to a cone over $M_t$,
\item the fibre $M_0$ is equal to $\SSS_{\epsilon}\cap\Sigma_0$.
\end{enumerate}
\end{lema}
\begin{proof}
The proof is an easy application (actually a simplification) of the technique of construction of cuts developed in~\cite{Bo}, Section 2.
\end{proof}
In this subsection we prove the following proposition:
\begin{prop}\label{existYcut}
Let $f_t$ be equisingular at the singular locus and equisingular at the normalisation.
Given any positive function $\eta:T\to \RR_+$ there exist a cut $\hat{H}$ for $\hat\Sigma$ such that $n(\mathrm{Int} \hat{H})\subset B(T,\eta)$.
\end{prop}
\begin{proof} $\quad$

We contruct first a cut for $\hat\Sigma$ over $0$ and then extend it over $T$ using the homeomorphism
$\alpha^{-1}$.

Choose $0<\epsilon_*<\eta(0)$ such that the sphere $\SSS_{\epsilon}$ is transverse (in the stratified sense) to $X_0$ for any $\epsilon\leq\epsilon_*$.

We can take $0<\epsilon_0<\epsilon_*$ such that
$\alpha^{-1}(T\times n_0^{-1}(B_{\epsilon_0}))$ is contained in $n^{-1}(B(T,\eta))$.

Using the technique of~\cite{Bo}, Section 3, we modify
$\SSS_{\epsilon_0}$ to get a neighbourhood $W$ of the origin in $\CC^3$ such that
\begin{itemize}
\item the hypersurface $\alpha^{-1}(T\times n_0^{-1}(\partial W))$ is contained in $n^{-1}(B(T,\eta))$,
\item we have the equality $\partial W\cap\Sigma_0=\SSS_{\epsilon_0}\cap \Sigma_0$,
\item the boundary $\partial W$ coincides with the normal bundle $E$ in a neighbourhood of $\Sigma_0\cap \partial W$ in $\partial W$. Precisely, we have the equality $$\partial W\cap E(\Sigma\setminus T,\mu)_0=E(\partial W\cap \Sigma_0,\mu).$$
\end{itemize}
The arguments needed here are much simpler than those in \cite{Bo} since we
have to ensure neither smoothness of $\partial W$ nor any cobordism condition.

The topological hypersurface
\begin{equation}\label{defL}L:=\alpha^{-1}(T\times n_0^{-1}(\partial W))\end{equation}
is a candidate for being a cut. 
It satisfies the properties $(i)$ and $(v)$ of Definition \ref{Ycut} and, over the parameter $0$, also the other ones.
The rest of the proof consists in modifying $\alpha$ so that the cut $L$, defined as in (\ref{defL}), satisfies all the required conditions over any $t\in T$.\\
For the rest of the proof we will need three auxiliary hypersurfaces of $\Sigma\setminus T$, that we will call $M_+$, $M_*$ and $M_-$, easily constructed applying Lemma~\ref{Singcut}, satisfying the following conditions:
\begin{itemize}
\item We have the inclussions
$$\mathrm{Int}_{\Sigma} M_-\subset \mathrm{Int}_{\Sigma} M_*\subset \mathrm{Int}_{\Sigma} M_+$$
$$\mathrm{Int}_{\hat\Sigma} \hat M_-\subset \mathrm{Int}_{\hat\Sigma} L \subset \mathrm{Int}_{\hat\Sigma}\hat M_+$$
\item The restriction $(M_*)_0$ coincides with $\SSS_{\epsilon_0}\cap \Sigma_0$,
\end{itemize}

Observe 
that $\hat M_-$, $\hat M_*$ and $\hat M_+$ are smooth topological hypersurfaces of $\hat\Sigma$ and that we have the equalities $\mathrm{Int}_{\hat\Sigma}\hat M_-=n^{-1}(\mathrm{Int}_\Sigma M_-)$, 
$\mathrm{Int}_{\hat\Sigma}\hat M_*=n^{-1}(\mathrm{Int}_{\Sigma} M_*)$ and
$\mathrm{Int}_{\hat\Sigma}\hat M_+=n^{-1}(\mathrm{Int}_{\Sigma} M_+)$.
 Observe also that the hypersurface $\hat M_*$ satisfies conditions $(i)-(iii)$ of Definition \ref{Ycut}.

For the rest of the proof we will assume that $\hat\Sigma$ is irreducible. The general case follows similarly. 

\noindent\textbf{Step 1}. We modify $\alpha$ so that $L$, defined as in (\ref{defL}), satisfies conditions $(i)-(iii)$ and $(v)$.
Define $P:=\overline{\mathrm{Int}\hat{M}_+\setminus \mathrm{Int}\hat{M}_-}$. The restriction $\hat\tau|_P:P\to T$ is a trivial
fibration with fibre a cylinder (in the general case, as many cylinders as irreducible components $\hat\Sigma$ has). This is true by condition $(iv)$ of Lemma~\ref{Singcut} and
the fact that the restriction $n|_{\hat\Sigma\setminus \hat T}$ is a covering projection.

Consider $Q:=n^{-1}(\pi|_{E(\Sigma\setminus T,\mu)}^{-1}(\overline{\mathrm{Int}M_+\setminus \mathrm{Int}M_-)})$ and the inherited bundle structure $\hat{\pi}|_Q:=Q\to P$. The fibre is a disc $D\subset \CC$ and the bundle is trivial because it is a orientable fibration over $P$, which has the homotopy type of $\SSS^1$. Since the bundle is trivial we identify $Q$ with $P\times D$, noting that the product structure is compatible with the fibration $\hat\tau$.

We show that there exists a homeomorphism
$$\Theta:Q\to Q$$ such that it preserves the fibres of $\hat\tau$, it is the identity over $Q_0$, it preserves $P$, it carries $L\cap P$ to $\hat{M}_*$ and it leaves $\partial Q$ fixed.

Suppose that this homeomorphism exists. We extend it with the identity outside $Q$ obtaining a homeomorphism $$\bar{\Theta}:\hat X\to \hat X.$$
Since $\bar\Theta$ transforms $L\cap P$ in $\hat M_*$, the topological hypersurface $\bar\Theta(L)$ satisfies Properties $(i)-(iii)$ and $(v)$ of Definition \ref{Ycut}.
Redefining $\alpha:=\bar\Theta\comp\alpha$ we can assume that $L$, defined as in (\ref{defL}), satisfies Properties $(i)-(iii)$ and $(v)$.

To construct $\Theta$, we construct first an isotopy $$\theta: P\times [0,1]\to P\times [0,1]$$ such that
it preserves fibres of $\hat{\tau}$,
the restriction $\theta_1$ is equal to $Id_P$,
the restriction $\theta_0$ carries $L\cap P$ to $\hat M_*$ and
it leaves $\partial P=\hat M_+\cup \hat M_-$ fixed.

Then $\Theta$ can be easily constructed using $\theta$ as follows.
For any $(p,z)\in P\times D=Q$ we define $$\Theta(p,z):=(\theta_{|z|}(p),z).$$

Let us construct $\theta$. First of all, consider a trivialisation of $\hat\tau|_P:P\to T$
$$\beta:P\to T\times P_0$$
and a product structure for $P_0$, 
$$\gamma:P_0\to \SSS^1\times[0,1]$$
that carries $(\hat M_+)_0$ to $\SSS^1\times\{0\}$, $(L\cap P)_0$ (which is equal to $(\hat M_*)_0$) to $\SSS^1\times\{\frac{1}{2}\}$ and $(\hat M_-)_0$ to $\SSS^1\times\{1\}$.

      \begin{center}
      \includegraphics[width=110mm]{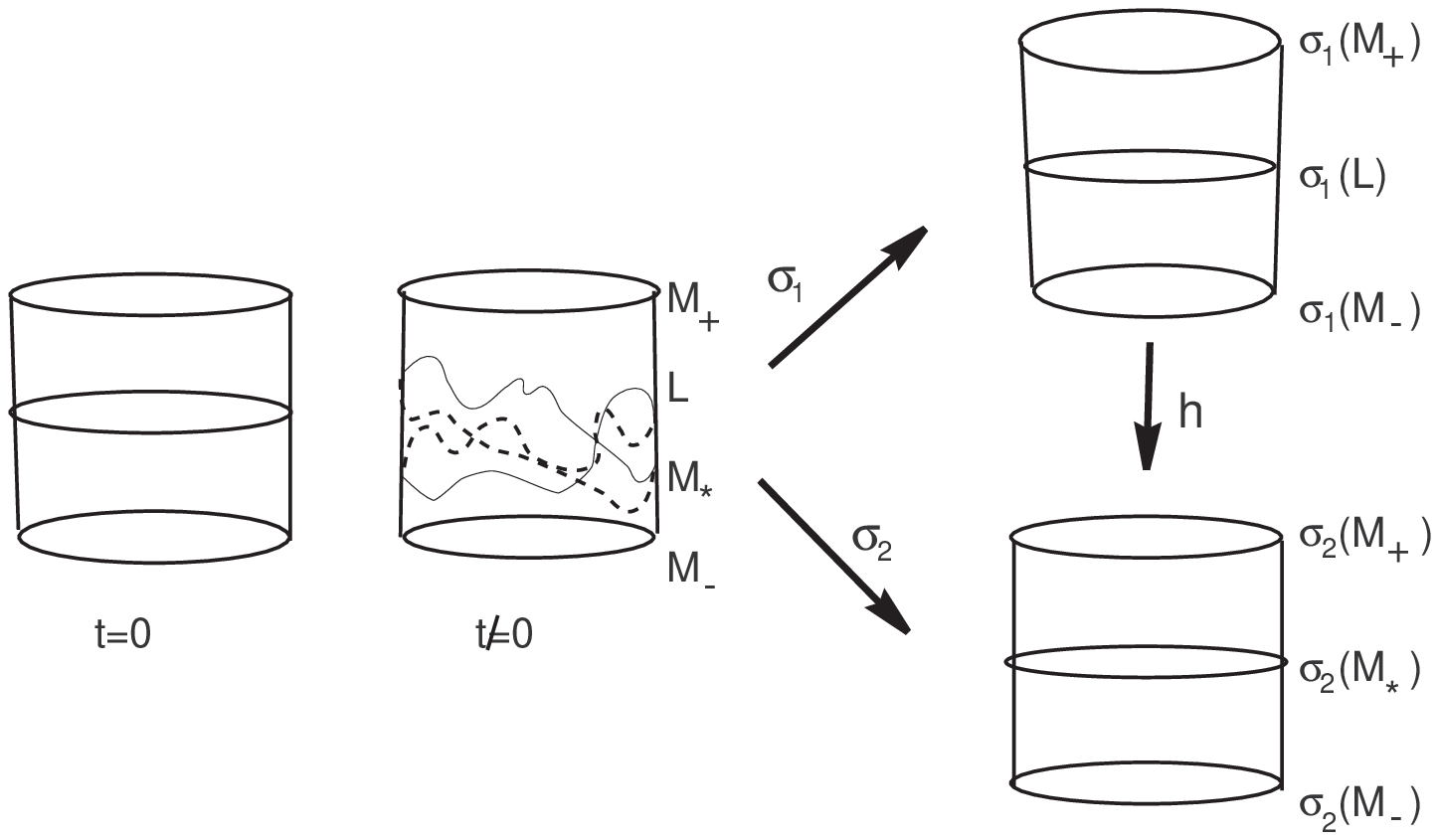}\label{conos}\\

      \end{center}

By construction, the hypersurface $L$ separates $P$ in two connected components, $P_{1,1}$ and $P_{1,2}$. Observe that both $\partial P_{1,1}$ and $\partial P_{1,2}$ have two connected components, $\hat M_+$ and $L$ for the first one, and  $L$ and $\hat M_-$ for the second one.
The restrictions $\hat\tau|_{P_{1,1}}$ and $\hat\tau|_{P_{1,2}}$ are trivial fibrations with fibre a cylinder so that we can take trivialisations
$\sigma_{1,1}:P_{1,1}\to T\times \SSS^1\times[0,\frac{1}{2}]$  which carries $M_+$ to $T\times \SSS^1\times\{0\}$ and $L$ to $T\times \SSS^1\times\{\frac{1}{2}\}$,
$\sigma_{1,2}:P_{1,2}\to T\times \SSS^1\times[\frac{1}{2},1]$ which carries $L$ to $T\times \SSS^1\times\{\frac{1}{2}\}$ and $\hat M_-$ to $T\times \SSS^1\times\{1\}$.

We can glue them along $T\times\SSS^1\times\{\frac{1}{2}\}$ with the mapping $\sigma_{1,2}\comp\sigma_{1,1}^{-1}$ to obtain a homeomorphism
$$\sigma_1:P\to T\times \SSS^1\times[0,1]=(T\times \SSS^1\times[0,\frac{1}{2}]){\bigcup}_{\sigma_{1,2}\comp\sigma_{1,1}^{-1}} (T\times \SSS^1\times[\frac{1}{2},1])$$
which carries $\hat M_+$ to $T\times \SSS^1\times\{0\}$, $L\cap P$ to $T\times \SSS^1\times\{\frac{1}{2}\}$ and $\hat M_-$ to $T\times \SSS^1\times\{1\}$.

The hypersurface $\hat M_*$ also separates $P$ in two connected components, $P_{2,1}$ and $P_{2,2}$. In this case, their frontiers have components $\hat M_+$ and $\hat M_*$, and $\hat M_*$ and $\hat M_-$ respectively. Analogously to the construction of $\sigma_1$ made for $L$ we can obtain a homeomorphism
$$\sigma_2:P\to T\times \SSS^1\times[0,1]$$
which carries $\hat M_+$ to $T\times \SSS^1\times\{0\}$, $\hat M_*$ to $T\times \SSS^1\times\{\frac{1}{2}\}$ and $\hat M_-$ to $T\times \SSS^1\times\{1\}$.

Consider the homeomorphism
$$h:=\beta\comp\sigma_2^{-1}\comp\sigma_1\comp\beta^{-1}:T\times P_0\to T\times P_0$$
which satisfies that the restriction $h_0$ is the identity, it carries $\beta(L\cap P)$ to $\beta(\hat M_*)$ and it preserves the fibres of the projection $p_1:T\times P_0\to T$.

Using Lemma \ref{glueing} we can modify $\sigma_2$ in collars of $M_+$, $M_-$ in $P$ such that $h$ also leaves $T\times (M_+)_0$ and $T\times (M_-)_0$ fixed.

Applying now Lemma \ref{glueprev}, there is an isotopy $$\phi:T\times P_0\times[0,1]\to T\times P_0\times [0,1]$$
which satisfies that $\phi_0=h$, $\phi_1=Id_{T\times P_0}$ and that $\phi_s$ preserves the fibres of the projection $p_1:T\times P_0\to T$
and leaves $T\times (M_+)_0$ and $T\times (M_-)_0$ fixed for any $s\in [0,1]$.

We finish defining $\theta$ as the composition
$$(\sigma_2^{-1},Id_{[0,1]})\comp (Id_T,\gamma,Id_{[0,1]})\comp \phi\comp(Id_T,\gamma^{-1},Id_{[0,1]})\comp(\sigma_1,Id_{[0,1]}).$$

\noindent\textbf{Step 2.} We modify again $\alpha$ in order to make $L$, defined as in (\ref{defL}), satisfy also condition $(iv)$ of Definition~\ref{Ycut}.

Take two tubular neighbourhoods $G_-\subset G_+$ of $L\cap \hat\Sigma_0$ in $\hat\Sigma_0$. Define $(K_-)_0:=\hat\pi^{-1}(G_-)\cap \hat E(\eta_-)$ and $(K_+)_0:=\hat\pi^{-1}(G_+)\cap \hat E(\eta_+)$ for positive functions $\eta_-,\eta_+:\Sigma\setminus T\to\RR^+$ such that $\eta_-<\eta_+<\mu$. Define
 $$K_+:=\alpha^{-1} (T\times (K_+)_0),$$
 $$K_-:=\hat\pi^{-1}(\alpha^{-1}(T\times G_-))\cap \hat E(\eta_-),$$ 
 $$L_*:=\hat\pi^{-1}(L\cap \hat\Sigma)\cap \hat E(\eta_-).$$

Notice that $(L_*)_0$ coincides with $L_0\cap \hat E(\eta_-)$. If $\eta_-$ is sufficiently small then we
have the inclussion $K_-\subset\dot K_+$.

      \begin{center}
      \includegraphics[width=75mm]{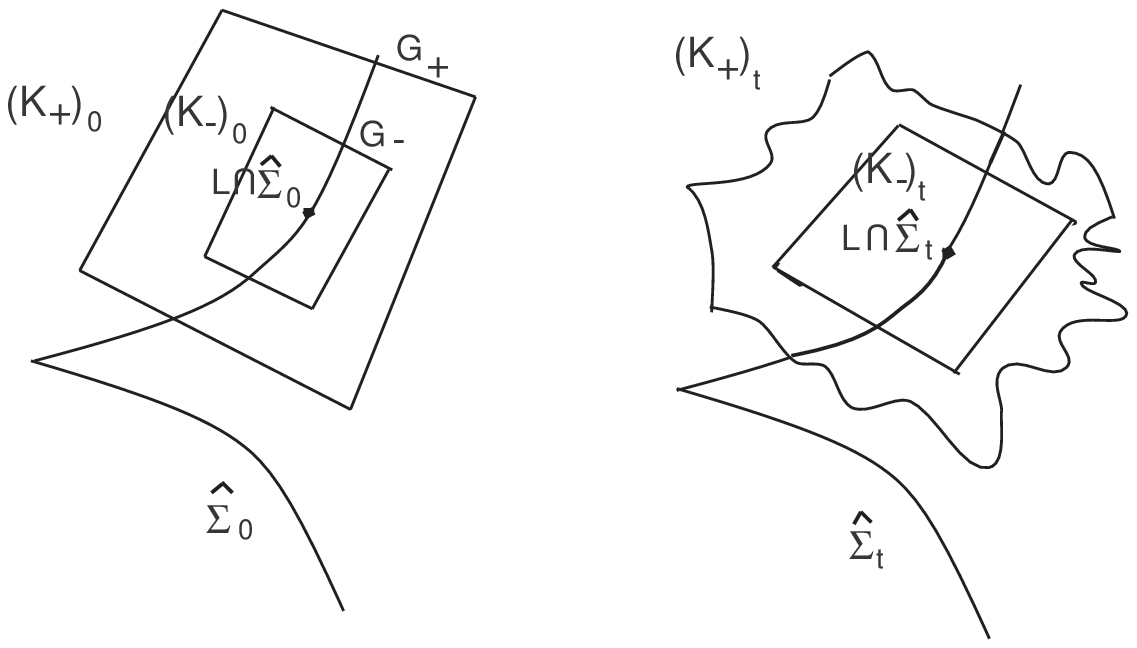}\label{corte_arrugado}\\

      \end{center}

By construction of $K_+$ we have a trivialisation  
 $$\alpha|_{\hat X\setminus K_+}:\hat X \setminus K_+\to T\times (\hat X\setminus K_+)_0.$$
of $\tau|_{\hat X\setminus K_+}$, and, as $\tau|_{K_+\setminus K_-}$ and  $\tau|_{K_-}$ are proper topological submersions, they admit trivialisations
 $$\kappa_1:(K_+\setminus K_-)\to T\times((K_+\setminus K_-)_0)$$
 $$\kappa_2:(K_-,L_*)\to T\times( (K_-)_0,(L\cap K_-)_0).$$ 

Taking collars of $\partial K_+$ and $\partial K_-$ in $K_+\setminus K_-$, we can apply Lemma \ref{glueing} and obtain a trivialisation
$$\kappa:\hat X\to T\times \hat X_0$$ such that $\kappa|_{\hat X\setminus K_+}=\alpha|_{\hat X\setminus K_+}$ and $\kappa|_{K_-}=\kappa_2$.

Redefining $\alpha:=\kappa$ and defining $\hat H$ as in (\ref{defL}), we have the cut requiered.

\end{proof}

\subsection{Equisingularity at the critical locus}\label{ssequiatcrit}

\begin{prop}
\label{singcrit}
If $f_t$ is equisingular at the singular locus and equisingular at the normalisation, then it is equisingular at the critical locus.
\end{prop}
\begin{proof}
We have to prove that the polar variety $\Gamma:=\overline{C\setminus\Sigma}$ is empty. It is clearly sufficient to assume that
$T$ is a disc. By the second property of equisingularity at the singular locus, 
Thus $\Gamma$ is a curve.
If $\Gamma$ does not meet $T$ we obtain the conclusion shrinking the representative $U$.
Suppose that $\Gamma$ meets $T$ at a point $t_0$. We may choose the base point $0$ to be equal to $t_0$. Denote by $\mu_0$ the
intersection multiplicity $I_{(0,O)}(\Gamma,\{0\}\times\CC^3)$.

Choose radii $\epsilon$, $\delta$ for the Milnor fibration of the function $f_0$ (in particular, $B_\epsilon\cap f_0^{-1}(D_\delta)$ is
contractible). The arguments given in the proof of Lemma 26 of~\cite{Bo} show that equisingularity
at the singular locus implies the existence of a small neighbourhood $T'$ of $0$ in $T$ such that $F^{-1}(t,s)$ is transversal to
$\SSS_\epsilon$ (in the stratified sense if $s=0$) for any $(t,s)\in T'\times D_\delta$. If $T'$ is small enough, for any $t\in T'\setminus\{0\}$, the function $f_t$ has
isolated critical points outside $f_t^{-1}(0)$ whose Milnor numbers add up $\mu_0$. Under these conditions it follows by standard techniques
that $X_t\cap B_\epsilon$ has the homotopy type of a bouquet of $\mu_0$ spheres of dimension $2$.

Using Proposition \ref{existYcut} we construct a cut $\hat{H}$ for $\hat \Sigma$ such that $\hat{H}_0$ is equal to $n_0^{-1}(\SSS_{\epsilon/2})$.
If $T'$ is small enough we have that $n(\mathrm{Int}\hat H)_t$ is contained in $\dot{B}_\epsilon$ for any $t\in T'$. We claim that the inclussion
$$n(\hat H)_t\hookrightarrow \overline{X_t\cap B_\epsilon\setminus n(\mathrm{Int}\hat H)_t}$$
is a homotopy equivalence.

If the claim is true, then $B_\epsilon\cap X_t$ has the same homotopy type as $X_t\cap n(\mathrm{Int}\hat H)_t$. 
For any $t\in T$, the pair $\mathrm{Int}\hat H\cap (\hat X_t,\hat \Sigma_t)$ is a
disjoint union of cones, whose vertices are the elements of $n_t^{-1}(t)$.
From the topological viewpoint the mapping 
$$n_t:\hat X_t\cap\mathrm{Int}\hat H\to n_t(\hat X_t\cap \mathrm{Int}\hat H)$$ is a quotienting of
$\hat X_t\cap \mathrm{Int}\hat H$ which can be factorised as the composition of two identifications
$$\hat X_t\cap \mathrm{Int}\hat H\stackrel{a}{\longrightarrow} Z\stackrel{b}{\longrightarrow}n_t(\hat X_t\cap \mathrm{Int}\hat H)$$ which have the following properties.
\begin{itemize}
\item The restriction
$a|_{\hat X_t\cap \mathrm{Int}\hat H\setminus \hat \Sigma_t}$ is injective. Let $\{\hat \Sigma_t^i\}_{i=1}^d$ be the irreducible components of $\hat \Sigma$. For any $i\leq d$ there is a cyclic
action on $\hat \Sigma_t^i$ of $\hat \Sigma_t$, being the vertex the only point with non-trivial isotropy
group such that the fibres of $a$ over any point of $a(\hat \Sigma_t^i)$ coincide with the orbits of the action.
\item The restriction $b|_{Z\setminus a(\hat \Sigma)}$ is injective. There is a certain subset of pairs $(i,j)$ with $i\neq j\leq d$, and a
homeomorphism $\gamma_{i,j}:a(\hat \Sigma_t^i)\to a(\hat \Sigma_t^j)$ which takes the image by $a$ of the vertex of $\hat \Sigma_t^i$ to the image by
$a$ of the vertex of $\hat \Sigma_t^j$ such that two distinct elemets $z_1,z_2\in Z$ are in the same fibre of $b$ if and only if
$\gamma_{i,j}(z_1)=z_2$ for a certain pair $(i,j)$.
\end{itemize}
Using the mapping cylinder of the identifications $a$ and $b$ and Mayer-Vietoris sequences we deduce that the second homology
group of $n_t(\hat X_t\cap \mathrm{Int}\hat H)$ coincides with the second homology group of $\hat X_t\cap \mathrm{Int}\hat H$,
which vanishes for being $\hat X_t\cap \mathrm{Int}\hat H$
a disjoin union of cones. This proves that $\mu_0$ is equal to $0$, and hence that $\Gamma$ is empty.

We only have to prove the claim. The claim is true for $t=0$ since $\hat H_0$ is equal to $n_0^{-1}(B_{\epsilon/2})$. Let us prove it for any $t\in T'$.
Define $R:=X\cap(T'\times B_\epsilon)\setminus n(\mathrm{Int}\hat H)$.
As the generic transversal Milnor number of $f_t$ at any point of $\Sigma\cap R$ only depends on the connected component
in which the point lies,
using the technique of Section 2 of~\cite{Bo} we can find a neighbourhood $V$ of $\Sigma\cap R$ in
$R$ such that $\partial_\pi V:=\cup_{x\in\Sigma\cap R}\partial V_x$ is a manifold which meets transversely $X_t$ for any $t$,
and a trivialisation
\[\Psi_1:V\to T'\times V_0\]
of $\tau|_{V}$. The restriction
$\tau|_{R\setminus\dot{V}}:R\setminus\dot{V}\to T'$
is a proper topological submersion. Thus, taking $T'$ contractible, it is a trivial fibration. Consider a trivialisation
\[\Psi_2:R\setminus\dot{V}\to T'\times (R\setminus\dot{V})_0.\]
By Lemma~\ref{glueing} we can modify $\Psi_2$ so that it glues with $\Psi_1$ to get a trivialisation of $\tau|_{R}$. Then $R_0$ is
homeomorphic to $R_t$ and the claim for $t\neq 0$ follows from the claim for $t=0$.
\end{proof}

\subsection{Existence of cuts for $\Sigma$}\label{ssXcut}

We prove that we can extend $n(\hat H)$, being $\hat H$ the cut 
constructed in Section \ref{ssYcut}, to obtain a cut $H$ for
$\Sigma$ such that $\mathrm{Int} H\subset B(T,\eta)$. We extend $n(\hat H)$
in two steps, first in a neighbourhood of $\Sigma$, and then out
of it.

\noindent\textbf{Step 1.} Define $$H_1:=E(n(\hat H\cap
\hat\Sigma),\mu).$$ Because of the property $(iv)$ of Definition
\ref{Ycut} of the cut $\hat{H}$, we can ensure the equality
\begin{equation} \label{H1} X\cap H_1=n(\hat H)\cap E(\Sigma,\mu)
\end{equation}

\noindent\textbf{Step 2.}
We can consider a complex smooth line subbundle $$\sigma:N\to X\setminus \Sigma$$ of the tangent bundle of $T\times \CC^3$ that verifies
\begin{itemize}
\item it is a normal bundle to $X$,
\item the fibres of $\sigma$ are contained in the fibres of $\tau$ and
\item for any point $x\in E(\Sigma, \mu)\cap n(\hat H)\setminus \Sigma$ we have the inclusion $N_x\subset E_{\pi(x)}$.
\end{itemize}

Since the restriction of $f$ to each fibre of $\sigma$ is a local
diffeomorphism at the origin (for being $f$ non-singular at $X\setminus\Sigma$ and $N$ a normal bundle),
we can consider the metric obtained,
in a neighbourhood of $X\setminus \Sigma$, as the pullback by $f$ of a euclidean
metric in $\CC$. We can take a positive function $\nu:X\setminus
\Sigma\to \RR_+$ such that $N(X\setminus \Sigma,\nu)$ is embedded
in $T\times \CC^3$ as a tubular neighbourhood of $X\setminus
\Sigma$ in $T\times \CC^3\setminus\Sigma$.

Consider the topological compact manifold with boundary $M:=\overline{n(\hat H)\setminus E(\Sigma,\mu/2)}$.
By compactness of $M$ 
 we can consider $\delta:=\min \{\nu (x) : x\in M\}$. By the definition of the metric in $N$, the fibres $f_t^{-1}(s)$ are transverse to $N(M,\delta)$ for any $(t,s)\in T\times D_\delta$. 
Define $$H_2:=N(M,\delta).$$

Thanks to Equality (\ref{H1}) and the third property of $\sigma$,
the subspace $H_1\cup H_2$ is a topological hypersuface of
$T\times \CC^3$. Define $$H:=H_1\cup H_2.$$ If we have taken $\mu$
and $\delta$ small enough, we can get $H\subset B(T,\eta)$. Choose
$\delta>0$ such that $f_t^{-1}(s)$ is transverse to $H$ (in the
stratified sense when $s=0$) for any $(t,s)\in T\times D_\delta$.

Property $(iv)$ of Definition \ref{Xcut} follows directly from the
properties of $\hat H$. Properties $(i)-(ii)$ are satisfied by
construction. Let us verify Property $(iii)$.

We have to prove that $\mathrm{Int}H(\delta)$ (defined as in Definition~\ref{interiores}) is a bounded 
neighbourhood of $T$, the closure $\overline{\mathrm{Int}H(\delta)}$ is a topological manifold and the 
restriction 
\[\tau|_{\overline{\mathrm{Int}H(\delta)}}:\overline{\mathrm{Int}H(\delta)}\to T\] 
a topological submersion. It is easy to see that it is enough to prove the statement locally over the 
base $T$. In the proof we will possibly have to shrink $\mu$ and $\delta$. As $T$ is compact such a 
shrinking only shall occur a finite number of times.

Take $t_0\in T$ and
$\epsilon>0$, $\delta>0$ such that $\SSS_\epsilon$ is transverse to
$f_{t_0}^{-1}(s)$ for any $s\in D_\delta$ and such that
$B_\epsilon\cap H_{t_0}=\emptyset$. Take a neighbourhood $T'$ of
$t_0$ in $T$ where these conditions are satisfied for any $t\in
T'$.

Consider 
$$ K:=f^{-1}(T'\times\partial D_\delta)\cup H. $$ It is clear
that $\overline{Int_{T'\times \CC^3} K}=\overline{\mathrm{Int} H(\delta)}$.

Hence it is enough to construct a topological
manifold $L$, included in $K$, which satisfies the following
properties:
\begin{enumerate}[(a)] 
\item The restriction $\tau|_L:L\to T'$ is a proper topological submersion whose fibre
is a compact connected orientable manifold $L_t$; thus, applying the Jordan-Brower Separation Theorem, 
for each $t\in T'$, the space $\CC^3\setminus L_t$ has two connected components. 
The bounded ones, that we call $\mathrm{Int} L_t$, are neighbourhoods of the
origin and $\partial \mathrm{Int}
L_t=L_t$.
\item The intersection $\overline {\mathrm{Int} L}\cap K$ coincides with
$L$; therefore $\overline {\mathrm{Int} L}=\overline {\mathrm{Int} K}$.
\end{enumerate}

Let us construct $L$. Define $$G_1:=({T'}\times B_\epsilon)\cap
f^{-1}(D_\delta),$$ $$L_1:=\partial G_1.$$ The restriction
$\tau|_{L_1}$ is a
proper topological submersion. 
Define $$G_2:=G_1\cup E(\mathrm{Int}_\Sigma(H\cap \Sigma),\mu)\cap f^{-1}(D_\delta),$$
$$L_2:=\partial G_2.$$ Because of the transversality of $\Sigma_t
$ and $\SSS_\epsilon$ for every $t\in T'$, we can choose $\mu$
small enough (reducing $\delta$ accordingly) so that
$L_2$ is a topological manifold and
$\tau|_{L_2}$ is a topological submersion. Moreover $X\setminus\Sigma$ is transversal to $L_2$.

Define $$G:= G_2\cup N_{T'}(M,\delta),$$ $$L:=\partial
G.$$

As $X\setminus\Sigma$ is transversal to $L_2$, if $\delta$ is small enough then $L$ is a topological 
manifold. By construction it satisfies $(a)$ and $(b)$.

We have proved the following lemma
\begin{lema}\label{existXcut}
Let $f_t$ be equisingular at the
singular locus and equisingular at the normalisation. Given any
positive function $\eta:T\to \RR_+$, there exists a cut $H$ of
amplitude $\delta$ such that $\mathrm{Int} H(\delta)$ is included in
$B(T,\eta)$. \end{lema}

\subsection{Proof of Theorem \ref{equitop}}
We restate and prove Theorem \ref{equitop} in the following way,

\begin{theo}\label{maintech}
Let $f_t$ be equisingular at the singular locus and equisingular at the normalization.
Consider a cut $H$ for $\Sigma$ of amplitude $\delta$.
There exists a trivialisation  of $\tau$ relative to $f$
$$\Theta:\mathrm{Int} H(\delta)\to T \times (\mathrm{Int} H)_0.$$
\end{theo}

\begin{proof}
Applying Lemma \ref{existXcut} we construct a family of cuts $H_i$ of amplitude $\delta_i$ (with $\delta_1>\delta_2>...$) such that
$$ \mathrm{Int} H(\delta)=\mathrm{Int} H_1(\delta_1)\supset ...\supset \mathrm{Int} H_i(\delta_i)\supset \mathrm{Int} H_{i+1}(\delta_{i+1})\supset ...$$
and $$\bigcap_{i\in \NN} \mathrm{Int} H_i(\delta_i)=T.$$

Define $$R_i:=\overline{\mathrm{Int} H_i(\delta_i)\setminus \mathrm{Int} H_{i+1}(\delta_{i+1})}.$$

Following \cite{Bo}, let $\partial_\tau \mathrm{Int} H_i(\delta_i)$ denote $\bigcup_{t\in T} \partial (\mathrm{Int} H_i(\delta_i))_t$.

\begin{lema}\label{trivR} 
Any trivialisation $$\Xi_i:\partial_\tau \mathrm{Int} H_i(\delta_i)\to T\times \partial(\mathrm{Int} H_i(\delta_i)_0)$$
of $\tau|_{\partial_\tau \mathrm{Int} H_i(\delta_i)}$ relative to $f$
can be extended to a trivialisation $$\Theta_i:R_i\to T\times (R_i)_0$$ of $\tau|_{R_i}$ relative to $f$.
\end{lema}
\begin{proof}[Proof of Lemma \ref{trivR}]
We construct the trivialisation $\Theta_i$ piecewise. Split $R_i$ in the following two spaces $$P:=\mathrm{Int} H_i(\delta_{i+1})\setminus \mathrm{Int} H_{i+1}(\delta_{i+1})$$  $$Q:=\mathrm{Int} H_i(\delta_{i})\setminus \mathrm{Int} H_{i}(\delta_{i+1}).$$

Consider $P$ as the union of two subsets, $A$ and $B:=\overline{P\setminus A}$ where $A$ is defined below.

The function $f$ restricted to any of the fibers of the complex vector bundle $$\pi|_{E(\Sigma\cap R_i,\mu)}:E(\Sigma\cap R_i,\mu)\to \Sigma\cap R_i$$ is a holomorphic function with an isolated singularity at the origin whose Milnor number is independent of the considered fibre.
By Section 2 of \cite{Bo} there exist a cut $S$ of amplitude $\delta_{i+1}$ (we would restrict $\delta_{i+1}$ if neccesary) satisfying the inclusion $$Y_{int}(S,\Sigma\cap R_i,\delta_{i+1})\subset E(\Sigma\cap R_i,\mu).$$
We warn the reader that the notion of cut used in this paper is different from the one used in \cite{Bo}; the cut $S$ is a cut in the sense of \cite{Bo} and following the notation of that paper, $Y_{int}(S,\Sigma\cap R_i,\delta_{i+1})$ stands for the union of the connected components of $E(\Sigma\cap R_i,\mu)\cap f^{-1}(D_{\delta_{i+1}})\setminus S$ meeting $\Sigma\cap R_i$. 

Define $A:=Y_{int}(S,\Sigma\cap R_i,\delta_{i+1})$. 
By Theorem 11 of \cite{Bo}, there is a trivialisation 
$$\Theta_A:A\to T\times A_0$$ of $\tau|_A$ relative to $f$.
For being $f_t$ equisingular at the critical locus, $$F|_B:(B,A\cap B)\to T\times D_{\delta_{i+1}}$$ is a topological proper submersion of pairs; thus, we can extend the trivialisation of $F|_{A\cap B}$ given by $\Theta_A$ to a trivialisation of  $F|_B$. 
This trivialisation induces a trivialisation of $\tau|_B$ relative to $f$
$$\Theta_B:B\to T\times B_0$$
whose restriction to $A\cap B$ coincides with $\Theta_A$.

Define $\Theta_P:P\to T\times P_0$ glueing $\Theta_A$ and $\Theta_B$.
Analogously, the restriction
$$F|_Q:Q \to T\times (\overline{D_{\delta_i}\setminus D_{\delta_{i+1}}})$$
is a locally trivial fibration and there exist a homeomorphism 
$$\beta:Q\to (Q\cap f^{-1}(\partial D_{\delta_i}))\times [0,1]$$
which preserves the fibres of $\tau$ and carries fibres of $f$ to fibres of $(f\comp p_1,p_2)$ being $p_i$ the projections from $(Q\cap f^{-1}(\partial D_{\delta_i}))\times [0,1]$ to the $i$th factor.

Consider the restrictions
\[\beta_0:Q_0\to Q_0\cap f^{-1}(\partial D_{\delta_i})\times [0,1]\]
\[\alpha_0:Q_0\cap f^{-1}(\partial D_{\delta_{i+1}})\to Q_0\cap f^{-1}(\partial D_{\delta_i})\times \{1\}\cong Q_0\cap f^{-1}(\partial D_{\delta_i}).\]

Using the trivialisation $\Xi_i$ of the frontier, we consider the composition
$$\gamma:= (\Xi_i|_{Q\cap f^{-1}(\partial D_{\delta_i})},Id_{[0,1]})\comp \beta: Q\to T\times (Q_0\cap f^{-1}(\partial D_{\delta_i}))\times [0,1].$$
Composing with $(Id_T,\beta_0^{-1})$ we obtain a trivialisation of $\tau|_Q$ relative to $f$ that extends $\Xi_i|_{f^{-1}(\partial D_{\delta_i})\cap\mathrm{Int}H_i(\delta_i)}$:
 $$\tilde \Xi_i:=(Id_T,\beta_0^{-1})\comp \gamma:Q\to T\times Q_0$$

We will modify $\tilde \Xi_i$ in several steps, so that its restriction to $P\cap Q$ 
coincides with $\Theta_P|_{P\cap Q}$.

We extend first $\Theta_P|_{P\cap Q}$ to $Q\cap f^{-1}(\partial D_{\delta_{i+1}})$.
For being $F|_{\mathrm{Int} H_{i+1}\cap f^{-1}(\partial D_{\delta_{i+1}})}$ a proper topological submersion, 
there is a trivialisation of $\tau|_{\mathrm{Int} H_{i+1}\cap f^{-1}(\partial D_{\delta_{i+1}})}$ relative to $f$ that we call
$$\Upsilon:\mathrm{Int} H_{i+1}\cap f^{-1}(\partial D_{\delta_{i+1}})\to T\times (\mathrm{Int} H_{i+1}\cap f^{-1}(\partial D_{\delta_{i+1}}))_0.$$

As it is clear that $(P\cap Q)_0\cap H_{i+1}$ has a collar $K$ in $(\mathrm{Int} H_{i+1})_0\cap f^{-1}(\partial D_{\delta_{i+1}})$ relative to $f$, applying Lemma \ref{glueing} to $\Upsilon$ and $\Theta_P|_{P\cap Q}$, we obtain a trivialisation of $\tau|_{\mathrm{Int} H_{i+1}\cap f^{-1}(\partial D_{\delta_{i+1}})}$ relative to $f$
$$\tilde \Upsilon:Q\cap f^{-1}(\partial D_{\delta_{i+1}})\to T\times (Q_0\cap f^{-1}(\partial D_{\delta_{i+1}})).$$
with $\tilde\Upsilon|_{P\cap Q}=\Theta_P|_{P\cap Q}$ and $\tilde\Upsilon|_{Q\cap f^{-1}(\partial D_{\delta_{i+1}})\setminus \Upsilon^{-1}(K)}=\Upsilon|_{Q\cap f^{-1}(\partial D_{\delta_{i+1}})\setminus \Upsilon^{-1}(K)}$.

In order that $\tilde\Xi_i$ can be glued with $\Upsilon$ it is enough that the self-homeomorphism of $T\times (Q_0\cap f^{-1}(\partial D_{\delta_{i+1}}))$, defined by 
$$h:=\tilde\Xi_i|_{Q\cap f^{-1}(\partial D_{\delta_{i+1}})}\comp\tilde\Upsilon^{-1}$$ 
is the identity. 

By construction we have that $h_0$ is the identity. By Lemma \ref{glueprev} there is an isotopy
$$\phi:T\times (Q_0\cap f^{-1}(\partial D_{\delta_{i+1}})\times [0,1]\to T\times (Q_0\cap f^{-1}(\partial D_{\delta_{i+1}})\times [0,1]$$
such that $\phi_0=Id$ and $\phi_1=h^{-1}$.
Substituting $\gamma$ by 
$$(Id_T,\alpha_0^{-1},Id_{[0,1]})\comp\phi\comp (Id_T,\alpha_0,Id_{[0,1]})\comp\gamma$$ 
in the definition of $\tilde\Xi_i$ we obtain that $\tilde\Xi$ extends $\tilde\Upsilon$.

Now, we can define $\Theta$ glueing $\Theta_P$ and $\tilde\Xi_i$.
As there is a collar of $H_{i}(\delta_i)_0$ in $(R_i)_0$ relative to $f$ and to $f^{-1}(\partial D_{\delta_{i}})\cap (R_i)_0$, we can apply Lemma~\ref{glueing} to modify $\Theta$ so that it extends $\Xi_i$.

\end{proof}
We construct $\Theta_1$ following the proof of the previous Lemma, obtaining $\Theta_Q$ directly from the locally trivial fibration $(\tau,f):Q\to T\times (D_{\delta_i}\setminus D_{\delta_{i+1}})$ where $T$ is contractible.  

Using the lemma we construct inductively a family of trivialisations $$\Theta_i:R_i\to T\times (R_i)_0$$
of $\tau|_{R_i}$ relative to $f$ which glue to a homeomorphism $$\Theta^*:\mathrm{Int} H(\delta)\setminus T\to T\times (\mathrm{Int} H(\delta)\setminus T)_0.$$

This homeomorphism can be extended to the desired trivialisation simply defining $\Theta|_T:=(\tau,O)$.

\end{proof}

\section{Equisingularity of parametrised surfaces}
Using the theory developed above, we are going to study topological $\calA$-equisingularity of holomorphic families of surface parametrisations.
By a parametrisation we understand a finite, generically one to one mapping from $\CC^2$ to its image in $\CC^3$.
A parametrisation necessarily coincides with the normalisations of its image.

Let $T$ be a contractible and compact complex manifold. A family of parametrisations over $T$ is a finite holomorphic germ 
$$\varphi:T\times (\CC^2,O)\to T\times (\CC^3,O),$$ 
such that $\varphi_t(x):=\varphi(t,x)$ is a parametrisation.

\begin{definition}
The family $\varphi_t$ is {\em $\calA$-topologically equisingular} if there are homeomorphism germs
\[\beta:T\times (\CC^2,O)\to T\times(\CC^2,O)\]
$$\gamma:T\times(\CC^3,O)\to T\times (\CC^3,O)$$
which preserve the fibres of the respective projections to $T$ and
such that we have the equality $$\gamma\comp\varphi\comp\beta=(Id_T,\varphi_0).$$
If in addition the restrictions $\beta|_{T\times \CC^2\setminus T}$, $\gamma|_{T\times\CC^3\setminus \Sigma}$ and $\gamma|_{\Sigma\setminus T}$ are smooth, we say that the family $\varphi_t$ is {\em $\calA$-piecewise smoothly equisingular}.
\end{definition}

\begin{lema}\label{C2}
Let $X\subset T\times\CC^3$ be a family of surfaces and $n:\hat X\to X$ be the normalisation.
If $\hat X_t=\CC^2$ for any $t\in T$, and $X_t$ is equisingular at the singular locus, then the family $\hat\Sigma_t$ is a reduced $\mu$-constant 
family of plane curves.
\end{lema}
\begin{proof}
By the second condition of equisingularity at the singular locus the space $X\setminus T$ is equisingular at $\Sigma$. Thus the restriction 
$n:T\times(\CC^2\setminus\{O\})\to X\setminus T$ is the simulaneous resolution of the singularities of $X\setminus T$. This implies that 
the restriction $$n:\hat\Sigma\setminus\hat T\to \Sigma\setminus T$$ is an unramified covering. Hence, the first
condition of equisingularity at the singular locus, implies that $\hat\Sigma\setminus\hat T$ is a topological trivial fibration over $T$. 

The conclusion of the Lemma is now easy to obtain.
\end{proof}

\begin{theo}\label{C2Le}
Let $\varphi_t$ be a family of parametrisations. If the family is topologically $\calA$-equisingular, then the Milnor number of $\hat\Sigma_t$ at the origin is constant. If L\^e's Conjecture is true the converse holds.

If for a single $t\in T$ each of the generic tranversal singularities of $\varphi_t(\CC^2)$ has at least one smooth branch, then the converse holds without L\^e's Conjecture.
\end{theo}

\begin{proof}
Let $f:T\times\CC^3\to \CC$ be a reduced equation of the image of $\varphi$. The topological $\calA$-equisingularity of $\varphi_t$ implies the embedded topological equisingularity of $f_t$. This implies that $f_t$ is abstractly equisingular in the strong sense.
By Proposition \ref{abstonorm} and Lemma \ref{consecnorm} $(vii)$, we have that the Milnor number of $\hat\Sigma_t$ is constant.

Assuming the constancy of the Milnor number of $\hat\Sigma_t$ and the fact that $\hat X_t$ is $\CC^2$ for any $t\in T$, we obtain that $f_t$ is equisingular at the normalisation. By Proposition \ref{equinormsing}, the family $f_t$ is equisingular at the singular locus. By Theorem \ref {maintech}, there is a family of homemorphisms $\gamma_t:(\CC^3,0)\to (\CC^3,0)$ such that $\gamma_t(X_0)=X_t$.
Using Lemma \ref{lift}, we can construct $\beta_t:(\CC^2,0)\to (\CC^2,0)$ such that $$\gamma_t\comp \varphi_t\comp \beta_t=\varphi_0.$$

\end{proof}

\begin{cor}
Let $\varphi_t$ be a family of parametrisations. If the generic tranversal singularities of $\varphi_t(\CC^2)$ are of Morse type, then the family is topologically $\calA$-equisingular if and only if the Milnor Number of $\hat\Sigma_t$ is independent of $t\in T$.
\end{cor}

\begin{cor}[Calleja-Bedregal, Houston, Ruas]
Let $\varphi_t$ be a family of finitely determined mapping germs.
The family is topologically $\calA$-equisingular if and only if the Milnor Number of $\hat\Sigma_t$ is independent of $t\in T$.
\end{cor}
\begin{proof}
In the finitely determined case, mapping germs are parametrisations and transversal singularities are of Morse type.
\end{proof}

\begin{theo}
Let $\varphi_t$ a family of parametrisations and $f_t$ the reduced equations of the $X_t$.
The constancy of the generic L\^e numbers of $f_t$ implies topological $\calA$-equisingularity of $\varphi_t$.
\end{theo}
\begin{proof}
The constancy of L\^e numbers implies equisingularity at the critical locus as it was proved in \cite{Bo} Theorem 48. By Lemma \ref{C2} the family $\hat\Sigma_t$ is $\mu$-constant and hence $\varphi_t$ is equisingular at the normalisation. Once we have equisingularity at the singular locus and equisingularity at the normalisation, we reproduce the proof of Theorem \ref{C2Le}.
\end{proof}

\section{The piecewise smooth equisingularity}

The tecnique of this article can be readily adapted to get the analogous results for the piecewise smooth case.

\begin{theo}\label{ps}
Let $f_t$ be equisingular at the singular locus and piecewise smooth equisingular at the normalisation,
then it is piecewise smooth $R$-equisingular.
\end{theo}

\begin{cor}\label{smoothC2Le}
Let $\varphi_t$ a family of parametrisations. If L\^e's Conjecture is true and the Milnor number of $\hat\Sigma_t$ at the origin is constant, then the family is piecewise smoothly $\calA$-equisingular.

If for a single $t\in T$ each of the generic tranversal singularities of $\varphi_t(\CC^2)$ has at least a smooth branch, then the conclusion holds without L\^e's Conjecture.
\end{cor}

\begin{proof} Since $\hat X_t=\CC^2$ for any $t\in T$, the constancy of Milnor number of $\hat\Sigma_t$, implies piecewise smooth equisingularity at the normalisation.
\end{proof}

Observe that L\^e's Conjecture implies that the notions of topological $\calA$-equisingularity and piecewise smooth $\calA$-equisingularity are equivalent since topological $\calA$-equisingularity implies the constancy of the Milnor number of $\hat\Sigma_t$.

\begin{cor}
Let $\varphi_t$ a family of parametrisations. If the generic tranversal singularities of $\varphi_t(\CC^2)$ are of Morse type, then the family is piecewise smoothly $\calA$-equisingular if and only if the Milnor Number of $\hat \Sigma_t$ is independent of $t\in T$.
\end{cor}

\begin{cor}
Let $\varphi_t:(\CC^2,0)\to(\CC^3,0)$ be a holomorphic family of finitely determined mapping germs.
The family is piecewise smoothly $\calA$-equisingular if and only if the Milnor Number of $\hat \Sigma_t$ is independent of $t\in T$.
\end{cor}

\end{document}